# Robust Inference of Risks of Large Portfolios

Jianqing Fan[*], Fang Han[†], Han Liu[‡], and Byron Vickers [§¶]

January 10, 2015

**Abstract**

We propose a bootstrap-based robust high-confidence level upper bound (Robust H-CLUB) for assessing the risks of large portfolios. The proposed approach exploits rank-based and quantile-based estimators, and can be viewed as a robust extension of the H-CLUB method (Fan et al., 2015). Such an extension allows us to handle possibly misspecified models and heavy-tailed data. Under mixing conditions, we analyze the proposed approach and demonstrate its advantage over the H-CLUB. We further provide thorough numerical results to back up the developed theory. We also apply the proposed method to analyze a stock market dataset.

## 1 Introduction

Let $\boldsymbol{R}_1, \ldots, \boldsymbol{R}_T$ be a stationary multivariate time series with $\boldsymbol{R}_t \in \mathbb{R}^d$ representing the asset returns at time $t$. Letting $\mathbf{w} \in \mathbb{R}^d$ be a portfolio allocation vector, we define the risk of $\mathbf{w}$ as

$$\text{Risk}(\mathbf{w}) := (\text{Var}(\mathbf{w}^\intercal \boldsymbol{R}_t))^{1/2} = (\mathbf{w}^\intercal \boldsymbol{\Sigma} \mathbf{w})^{1/2},$$

[*]Department of Operations Research and Financial Engineering, Princeton University, Princeton, NJ 08544, USA; e-mail: jqfan@princeton.edu. His research is supported by NSF grant DMS-1406266 and NIH grant R01GM100474-04.

[†]Department of Biostatistics, Johns Hopkins University, Baltimore, MD 21205, USA; e-mail: fhan@jhu.edu. His research is supported by a Google fellowship.

[‡]Department of Operations Research and Financial Engineering, Princeton University, Princeton, NJ 08544, USA; e-mail: hanliu@princeton.edu. His research is supported by NSF CAREER Award DMS1454377, NSF IIS1408910, NSF IIS1332109, NIH R01MH102339, NIH R01GM083084, and NIH R01HG06841.

[§]Department of Operations Research and Financial Engineering, Princeton University, Princeton, NJ 08544, USA; e-mail: bvickers@princeton.edu. His research was supported by NIH 2R01-GM072611-10.

[¶]We thank Huitong Qiu for discussions.

where $\mathbf{\Sigma}$ denotes the unknown volatility (or covariance) matrix of $\boldsymbol{R}_t$. i.e.,

$$\mathbf{\Sigma} := \mathbb{E}\big[(\boldsymbol{R}_t - \mathbb{E}\boldsymbol{R}_t)(\boldsymbol{R}_t - \mathbb{E}\boldsymbol{R}_t)^{\mathsf{T}}\big].$$

Assessing the risk of a portfolio includes two steps: First, we need a covariance matrix estimator $\widehat{\mathbf{\Sigma}}_{\text{est}}$; Secondly, we construct a confidence interval for $\mathbf{w}^{\mathsf{T}}\mathbf{\Sigma}\mathbf{w}$ based on $\widehat{\mathbf{\Sigma}}_{\text{est}}$.

Assessing the risk Risk($\mathbf{w}$) is challenging when $d$ is large. For example, given a pool of 2,000 candidate assets, the volatility matrix $\mathbf{\Sigma}$ involves more than 2 million parameters. However, for daily returns data, the sample size is in general no larger than 500 over one and a half years. This is a typical "small $n$, large $d$" problem which leads to the accumulation of estimation errors (Jagannathan and Ma, 2003; Pesaran and Zaffaroni, 2008; Fan et al., 2012). To handle the curse of dimensionality, more structural regularization is imposed in estimating $\mathbf{\Sigma}$. For example, Fan et al. (2008) and Fan et al. (2013) impose the factor model structure on the covariance matrix. The assumed factor structure reduces the effective number of parameters that have to be estimated. In addition, Ledoit and Wolf (2003) propose a shrinkage estimator of $\mathbf{\Sigma}$. Moreover, Barndorff-Nielsen (2002), Zhang et al. (2005), and Fan et al. (2012) consider estimating $\mathbf{\Sigma}$ based on high-frequency data. Other literature includes Chang and Tsay (2010), Gómez and Gallón (2011), Lai et al. (2011), Fan et al. (2011), Bai and Liao (2012), and Fryzlewicz (2013).

However, most of these papers focus on risk estimation instead of uncertainty assessment. To construct a confidence interval for $\mathbf{w}^{\mathsf{T}}\mathbf{\Sigma}\mathbf{w}$, Fan et al. (2012) propose to use $\|\mathbf{w}\|_1^2\|\widehat{\mathbf{\Sigma}}_{\text{est}} - \mathbf{\Sigma}\|_{\max}$[1] as an upper bound of $|\mathbf{w}^{\mathsf{T}}(\widehat{\mathbf{\Sigma}}_{\text{est}} - \mathbf{\Sigma})\mathbf{w}|$. However, this bound depends on the unknown $\mathbf{\Sigma}$ and has proven to be overly conservative in numerical studies. To handle this problem, Fan et al. (2015) further exploit several sample covariance based estimators $\widehat{\mathbf{\Sigma}}_{\text{est}}$ of $\mathbf{\Sigma}$ and propose a high-confidence level upper bound (H-CLUB) of $|\mathbf{w}^{\mathsf{T}}(\widehat{\mathbf{\Sigma}}_{\text{est}} - \mathbf{\Sigma})\mathbf{w}|$: For a given confidence level $1-\gamma$, under certain moment and dependence assumptions on the time series, the derived H-CLUB proves to dominate $|\mathbf{w}^{\mathsf{T}}(\widehat{\mathbf{\Sigma}}_{\text{est}} - \mathbf{\Sigma})\mathbf{w}|$ with probability approximating $1-\gamma$ as both $T$ and $d$ increase to infinity.

This paper proposes new methods for uncertainty assessment of risks of large portfolios for high dimensional heavy-tailed data. In particular, we derive confidence intervals for $\mathbf{w}^{\mathsf{T}}\mathbf{\Sigma}\mathbf{w}$ when the asset returns $\boldsymbol{R}_1, \ldots, \boldsymbol{R}_T$ are elliptically distributed. This setting has been commonly adopted in financial econometrics (Cont, 2001). To handle heavy-tailed data, we propose a new risk uncertainty assessment method named robust high-confidence level upper bound (Robust H-CLUB). The Robust H-CLUB exploits a new block-bootstrap-based approach for uncertainty assessment of Risk($\mathbf{w}$). More specifically, we decompose the problem of assessing the risk $\mathbf{w}^{\mathsf{T}}\mathbf{\Sigma}\mathbf{w}$ into two parts: (i) We propose a robust estimator $\widehat{\mathbf{\Sigma}}_{\text{est}}$ of $\mathbf{\Sigma}$; (ii) We derive the variance of $\mathbf{w}^\mathsf{T}(\widehat{\mathbf{\Sigma}}_{\rm est} - \mathbf{\Sigma})\mathbf{w}$. For estimating $\mathbf{\Sigma}$, we exploit rank-based Kendall's tau estimators and quantile-based median absolute deviation estimators. For estimating the variance of $\mathbf{w}^\mathsf{T}(\widehat{\mathbf{\Sigma}}_{\rm est} - \mathbf{\Sigma})\mathbf{w}$, we employ the circular block bootstrap method (Politis and Romano, 1992).

[1] We will provide the definitions of the vector $\ell_1$ norm ($\|\cdot\|_1$) and matrix $\ell_{\max}$ norm $\|\cdot\|_{\max}$ later.

Theoretically, when $T, d \to \infty$ and $d$ is possibly much larger than $T$, we develop an inferential theory of the robust risk estimators. In particular, we show that $\sqrt{T}\mathbf{w}^\mathsf{T}(\widehat{\mathbf{\Sigma}}_{\rm est} - \mathbf{\Sigma})\mathbf{w}$ is asymptotically normal with variance $\sigma^2$, and the block-bootstrap-based estimator $\widehat{\sigma}^2_{\rm est}$ of $\sigma^2$ is consistent. The theory holds even when $d$ is nearly exponentially larger than $T$. Moreover, it holds under any elliptical model. Thus we no longer need strong moment conditions (e.g., exponentially decaying rate on the tails of distributions) on the asset returns.

## 1.1 Other Related Work

There is a vast literature on estimating large sparse/factor-based covariance matrices. Under the assumption that data points are mutually independent, many sample covariance based regularization methods, including banding (Bickel and Levina, 2008b), tapering (Cai et al., 2010), thresholding (Bickel and Levina, 2008a; Cai and Zhou, 2012), and factor structures (Fan et al., 2008; Agarwal et al., 2012; Hsu et al., 2011), have been proposed. They are further applied to study stationary time series data under vector autoregressive dependence (Loh and Wainwright, 2012; Han and Liu, 2013c), mixing conditions (Pan and Yao, 2008; Fan et al., 2011, 2013; Han and Liu, 2013b), and physical dependence (Xiao and Wu, 2012; Chen et al., 2013).

This paper is also related to the literature on estimating large correlation/covariance matrix under the misspecified or heavy-tailed model. For example, Han and Liu (2014b), Han and Liu (2013a), Wegkamp and Zhao (2013), Mitra and Zhang (2014), and Fan et al. (2014) exploit the rank statistics, while Qiu et al. (2014) focus on quantile statistics. None of these works study the risk inference problem as in our paper.

## 1.2 Notation

Let $\mathbf{v} = (v_1, \ldots, v_d)^\mathsf{T}$ be a $d$ dimensional real vector and $\mathbf{M} = [M_{jk}]$ be a $d$ by $d$ real matrix. For $0 < q < \infty$, let the vector $\ell_q$ norm be $\|\mathbf{v}\|_q := (\sum_{j=1}^d |v_j|^q)^{1/q}$ and the vector $\ell_\infty$ norm be $\|\mathbf{v}\|_\infty := \max_{j=1}^d |v_j|$. For two subsets $I, J \in \{1, \ldots, d\}$, we denote $\mathbf{v}_I$ and $\mathbf{M}_{I,J}$ as the sub-vector of $\mathbf{v}$ with entries indexed by $I$ and sub-matrix of $\mathbf{M}$ with rows and columns indexed by $I$ and $J$. We denote the matrix $\ell_{\max}$ norm of $\mathbf{M}$ as $\|\mathbf{M}\|_{\max} := \max_{jk} |M_{jk}|$. Letting $\mathbf{N} = [N_{jk}] \in \mathbb{R}^{d\times d}$ be another $d$ by $d$ real matrix, we denote by $\mathbf{M} \circ \mathbf{N} = [M_{jk}N_{jk}]$ the Hadamard product between $\mathbf{M}$ and $\mathbf{N}$. Letting $f : \mathbb{R} \to \mathbb{R}$ be a real function, we denote by $f(\mathbf{M}) = [f(M_{jk})]$ the matrix with $f(M_{jk})$ as its $(j, k)$ entry. We write $\mathbf{M} =$ $\text{diag}(\mathbf{M}_1, \ldots, \mathbf{M}_k)$ if $\mathbf{M}$ is block diagonal with diagonal matrices $\mathbf{M}_1, \ldots, \mathbf{M}_k$. For random vectors $\boldsymbol{X}, \boldsymbol{Y} \in \mathbb{R}^d$, we write $\boldsymbol{X} \overset{\mathsf{d}}{=} \boldsymbol{Y}$ if $\boldsymbol{X}$ and $\boldsymbol{Y}$ are identically distributed. Throughout the paper, we use $c, c_1, c_2, \ldots$, and $C, C_1, C_2, \ldots$ to represent generic absolute positive constants, for which the actual values may change at from one line to another. For any real positive sequences $\{a_n\}$ and $\{b_n\}$, we write $a_n \gtrsim b_n$ if we have $a_n \geq cb_n$ for some absolute constant $c$ and all large enough $n$. We write $a_n \lesssim b_n$ if we have $b_n \gtrsim a_n$, and $a_n \asymp b_n$ if $a_n \lesssim b_n$ and $a_n \gtrsim b_n$. For $a \in \mathbb{R}$, we define $\lceil a \rceil$ and $\lfloor a \rfloor$ to be the smallest integer larger than $a$ and the largest integer smaller than $a$ respectively.

### 1.3 Paper Organization

The rest of this paper is organized as follows. Section 2 introduces the Robust H-CLUB estimator for assessing the uncertainty of the portfolio risk. We consider three settings: (i) The marginal variances of the returns are known; (ii) The marginal variances are unknown, but with additional information for helping determine the values; (iii) The marginal variances are unknown and there is no additional information available. Section 3 presents the inferential theory for the risk estimators and justifies the use of Robust H-CLUB. Sections 4 and 5 present synthetic and real data analyses to back up the developed theory. Section 6 summarizes the results and discusses future work. Section 7 presents all the proofs.

## 2 Robust H-CLUB

This section introduces the Robust H-CLUB method. We consider a multivariate time series of asset returns $\boldsymbol{R}_1, \ldots, \boldsymbol{R}_T$ with $\boldsymbol{R}_t = (R_{t1}, \ldots, R_{td})^\mathsf{T} \in \mathbb{R}^d$ for $t = 1, \ldots, T$. Let $\boldsymbol{\Sigma} := \text{Cov}(\boldsymbol{R}_t)$ be the covariance matrix and $\mathbf{D} \in \mathbb{R}^{d \times d}$ be a diagonal matrix with diagonals $\boldsymbol{\Sigma}_{11}^{1/2}, \ldots, \boldsymbol{\Sigma}_{dd}^{1/2}$. It is easy to derive $\boldsymbol{\Sigma} = \mathbf{D}\boldsymbol{\Sigma}^0\mathbf{D}$, where $\boldsymbol{\Sigma}^0$ is the correlation matrix of $\boldsymbol{R}_t$. For a given portfolio allocation vector $\mathbf{w} \in \mathbb{R}^d$, we aim to construct a confidence interval for $\mathbf{w}^\mathsf{T}\boldsymbol{\Sigma}\mathbf{w}$. Throughout this section, our interest is on analyzing heavy-tailed returns, which are common in financial applications.

We exploit the elliptical distribution family to model heavy-tailed data. The elliptical distribution is routinely used in modeling financial data (Owen and Rabinovitch, 1983; Hamada and Valdez, 2004; Frahm and Jaekel, 2007). More specifically, a random vector $\boldsymbol{Z} \in \mathbb{R}^d$ follows an elliptical distribution with mean $\boldsymbol{\mu} \in \mathbb{R}^d$ and positive definite covariance matrix $\boldsymbol{\Sigma} \in \mathbb{R}^{d \times d}$ if

$$\boldsymbol{Z} \overset{\mathsf{d}}{=} \boldsymbol{\mu} + \xi \mathbf{A}\boldsymbol{U},$$

where $\mathbf{A} \in \mathbb{R}^{d \times d}$ satisfies $\mathbf{A}\mathbf{A}^\mathsf{T} = \boldsymbol{\Sigma}$, $\boldsymbol{U} \in \mathbb{R}^d$ is uniformly distributed on the $d$-dimensional sphere $\mathbb{S}^{d-1}$, and $\xi$ is an unspecified nonnegative random variable independent of $\boldsymbol{U}$ satisfying

$\mathbb{E}\xi^2 = d$. We impose the following stationary assumption on $\{\boldsymbol{R}_t\}_{t=1}^T$:

- **(A0).** $\boldsymbol{R}_1, \ldots, \boldsymbol{R}_T$ are continuous and identically distributed as an elliptical random vector $\boldsymbol{R}$ with covariance and correlation matrices $\boldsymbol{\Sigma}$ and $\boldsymbol{\Sigma}^0$.

For parameter estimation, we define the rank-based Kendall's tau correlation coefficient and quantile-based median absolute deviation estimators. In detail, given $\boldsymbol{R}_1, \ldots, \boldsymbol{R}_T$, the sample and population Kendall's tau matrices $\widehat{\mathbf{T}} = [\widehat{\tau}_{jk}]$ and $\mathbf{T} = [\tau_{jk}]$ are defined as

$$\widehat{\tau}_{jk} := \frac{2}{T(T-1)} \sum_{t<t'} \text{sign}(R_{tj} - R_{t'j}) \text{sign}(R_{tk} - R_{t'k}),$$
$$\tau_{jk} := \mathbb{E}\text{sign}(R_j - \widetilde{R}_j)\text{sign}(R_k - \widetilde{R}_k), \tag{2.1}$$

where $\mathbf{R} = (R_1, \ldots, R_d)^\intercal$ and $\widetilde{\mathbf{R}} = (\widetilde{R}_1, \ldots, \widetilde{R}_d)^\intercal$ are two independent copies of $\mathbf{R}_1$. Under the elliptical model, the Kendall's tau matrix $\mathbf{T}$ and correlation matrix $\boldsymbol{\Sigma}^0$ satisfy (Lindskog et al., 2003):

$$\boldsymbol{\Sigma}^0_{jk} = \sin\Big(\frac{\pi}{2}\tau_{jk}\Big). \tag{2.2}$$

Next, we define the quantile-based median absolute deviation estimator of the scale parameter. We start with some extra notation. Let $X \in \mathbb{R}$ be a random variable and $\{X_1, \ldots, X_T\}$ be $T$ realizations of $X$. For any $q \in [0, 1]$, we define the population and sample $q$-quantiles as

$$Q(X; q) := \inf\big\{x : \mathbb{P}(X \le x) \ge q\big\},$$
$$\widehat{Q}(\{X_t\}; q) := X^{(k)}, \text{ where } k = \min\Big\{t : \frac{t}{T} \ge q\Big\}. \tag{2.3}$$

Here $X^{(1)} \le X^{(2)} \le \cdots \le X^{(T)}$ are the ordered sequence of $X_1, \ldots, X_T$[2]. We then define the population and sample median absolute deviations for $\{X_1, \ldots, X_T\}$ as the population and sample medians of absolute values of the centered data. The formal definitions are as follows:

$$\sigma_{\text{M}}(X) := Q\bigg(\Big\{\Big|X - Q\Big(X; \frac{1}{2}\Big)\Big|\Big\}; \frac{1}{2}\bigg),$$
$$\widehat{\sigma}_{\text{M}}(\{X_t\}_{t=1}^T) := \widehat{Q}\bigg(\Big\{\Big|X_t - \widehat{Q}\Big(\{X_t\}_{t=1}^T; \frac{1}{2}\Big)\Big|\Big\}_{t=1}^T; \frac{1}{2}\bigg). \tag{2.4}$$

They are robust alternatives to the population and sample standard deviations. In particular, for an elliptically distributed random vector $\boldsymbol{R} = (R_1, \ldots, R_d)^\intercal$, Han et al. (2014) prove that

$$\frac{\sigma_M(R_1)}{\text{sd}(R_1)} = \frac{\sigma_M(R_2)}{\text{sd}(R_2)} = \cdots = \frac{\sigma_M(R_d)}{\text{sd}(R_d)}, \tag{2.5}$$

[2] Let $F$ and $f$ be the distribution function and density function of $X$. We will use $Q(X; q), Q(F; q)$, and $Q(f; q)$ exchangeably.

where for arbitrary random variable $X$, $\mathrm{sd}(X)$ represents the standard deviation of $X$.

Under the elliptical model and using the rank- and quantile-based estimators, we propose three robust approaches to construct the confidence interval of $\mathbf{w}^\top\mathbf{\Sigma}\mathbf{w}$. Formally speaking, for each proposed robust covariance matrix estimator $\widehat{\mathbf{\Sigma}}_{\text{est}}$ and any given $\gamma > 0$, we aim to find a $\widehat{U}_{\text{est}}(\gamma)$ such that

$$\mathbb{P}\Big(\mathbf{w}^\top\mathbf{\Sigma}\mathbf{w} \in \big[\mathbf{w}^\top\widehat{\mathbf{\Sigma}}_{\text{est}}\mathbf{w} - \widehat{U}_{\text{est}}(\gamma), \mathbf{w}^\top\widehat{\mathbf{\Sigma}}_{\text{est}}\mathbf{w} + \widehat{U}_{\text{est}}(\gamma)\big]\Big) \to 1-\gamma,$$

as $T, d \to \infty$. The proposed approaches correspond to three scenarios where $\mathbf{D}$ has different structures.

Of note, a main strategy throughout the proposed three methods is to separately estimate the marginal standard deviations and bivariate correlation coefficients. In this paper, we focus on measuring the uncertainty introduced in estimating the correlation coefficients, while assuming that the uncertainty introduced in estimating marginal standard deviations is negligible[3]. For measuring the uncertainty in correlation coefficients estimation, we employ a circular block bootstrap method.

In detail, suppose that we derive robust marginal standard deviation estimator $\widehat{\mathbf{D}}_{\text{est}}$ of $\mathbf{D}$. We further derive the correlation matrix estimator $\widehat{\mathbf{\Sigma}}^0_{\text{est}}$ of $\mathbf{\Sigma}^0$ based on a $d$-dimensional multivariate time series $\boldsymbol{X}_1, \ldots, \boldsymbol{X}_T$. For any given portfolio allocation vector $\mathbf{w}$, we propose to estimate $\mathbf{w}^\top\mathbf{\Sigma}\mathbf{w}$ by

$$\widehat{\mathrm{Risk}}(\mathbf{w}) := \mathbf{w}^\top\widehat{\mathbf{\Sigma}}_{\text{est}}\mathbf{w}, \quad \text{where } \widehat{\mathbf{\Sigma}}_{\text{est}} := \widehat{\mathbf{D}}_{\text{est}}\widehat{\mathbf{\Sigma}}^0_{\text{est}}\widehat{\mathbf{D}}_{\text{est}}. \tag{2.6}$$

To estimate the asymptotic variance of the estimator $\mathbf{w}^\top\widehat{\mathbf{\Sigma}}_{\text{est}}\mathbf{w}$, we adopt a circular block bootstrap procedure introduced in Politis and Romano (1992). First, we extend the sample $\boldsymbol{X}_1, \ldots, \boldsymbol{X}_T$ periodically by concatenating $\boldsymbol{X}_{i+T} = \boldsymbol{X}_i$ for $i \geq 1$. We then randomly select a block of $l = l_T \asymp T^{1-\epsilon_0}$ consecutive observations from the extended sample for some absolute constant $\epsilon_0 < 1$ (e.g., we can pick $\epsilon_0$ to be 0.9). As the financial time series admits weakly dependence structure, the choice of block size $l$ is not very important. We repeat this process $b = \lfloor T/l \rfloor$ times independently to obtain a sample $\boldsymbol{X}^*_1, \ldots, \boldsymbol{X}^*_T$, so that for each $k = 0, \ldots, b-1$,

$$\mathbb{P}^*\big(\boldsymbol{X}^*_{kl+1} = \boldsymbol{X}_j, \ldots, \boldsymbol{X}^*_{(k+1)l} = \boldsymbol{X}_{j+l-1}\big) = 1/T, \text{ for } j = 1, \ldots, T,$$

where $\mathbb{P}^*$ is the resampling distribution conditional on $\boldsymbol{X}_1, \ldots, \boldsymbol{X}_T$. Based on each resampled time series $\boldsymbol{X}^*_1, \ldots, \boldsymbol{X}^*_T$, we calculate the correlation matrix estimator $\widehat{\mathbf{\Sigma}}^{0*}_{\text{est}}$. Let $\widehat{\mathbf{\Sigma}}^*_{\text{est}} := \widehat{\mathbf{D}}_{\text{est}}\widehat{\mathbf{\Sigma}}^{0*}_{\text{est}}\widehat{\mathbf{D}}_{\text{est}}$ be the estimator of $\mathbf{\Sigma}$ based on the resampled data and $\mathrm{Var}^*(\cdot)$ be the variance operator of the probability mass function $\mathbb{P}^*$. We estimate the asymptotic variance of $\mathbf{w}^\top \widehat{\mathbf{\Sigma}}_{\text{est}} \mathbf{w}$ by

$$\widehat{\sigma}^2_{\text{est}} := \text{Var}^*(\sqrt{T}\mathbf{w}^\top \widehat{\mathbf{\Sigma}}^*_{\text{est}} \mathbf{w}).$$

[3]This is mainly for the purpose of constructing the bootstrap-based inferential theory.

## 2.1 Known Marginal Volatilities

In this section we consider the setting where the marginal standard deviations of $\boldsymbol{R}_t$, encoded in $\mathbf{D}$, are known. While this is an ideal assumption, a practical implementation is to fit a parametric model such as the GARCH(1,1) model introduced in Bollerslev (1986) to each individual return time series. Such estimates are much more accurate than the nonparametric ones and can be ideally treated as known.

When $\mathbf{D}$ is known, estimating $\mathbf{w}^\top \mathbf{\Sigma} \mathbf{w}$ reduces to estimating the correlation matrix $\mathbf{\Sigma}^0$. Using (2.2), under the elliptical model, we focus on the covariance matrix estimator $\widehat{\mathbf{\Sigma}}$ with $\widehat{\mathbf{\Sigma}} := \mathbf{D} \sin\Big(\pi \widehat{\mathbf{T}}/2\Big)\mathbf{D}$. We then estimate $\mathbf{w}^\top \mathbf{\Sigma} \mathbf{w}$ via replacing $\widehat{\mathbf{\Sigma}}_{\text{est}}$ by $\widehat{\mathbf{\Sigma}}$ in (2.6). Let $\widehat{\sigma}^2$ be an estimator of the asymptotic variance $\sigma^2$ of $\mathbf{w}^\top \widehat{\mathbf{\Sigma}} \mathbf{w}$. We calculate $\widehat{\sigma}^2$ based on the circular block bootstrap method introduced earlier. Let $\Phi(\cdot)$ be the cumulative distribution function of a standard Gaussian random variable. For any given confidence level $1-\gamma \in (0,1)$, we define the Robust H-CLUB estimator $\widehat{U}(\gamma)$ as

$$\widehat{U}(\gamma) := \Phi^{-1}(1-\gamma/2)\sqrt{\widehat{\sigma}^2/T}. \tag{2.7}$$

The corresponding confidence interval for the risk is

$$\big[\mathbf{w}^\top \widehat{\mathbf{\Sigma}} \mathbf{w} - \widehat{U}(\gamma), \mathbf{w}^\top \widehat{\mathbf{\Sigma}} \mathbf{w} + \widehat{U}(\gamma)\big]. \tag{2.8}$$

In Section 3 we will show that, under mild conditions,

$$\widehat{\sigma}^2 = \sigma^2(1+o_P(1)) \quad \text{and} \quad \mathbb{P}\Big\{|\mathbf{w}^\top(\widehat{\mathbf{\Sigma}}-\mathbf{\Sigma})\mathbf{w}| \le \widehat{U}_\tau(\gamma)\Big\} \to 1-\gamma,$$

as $T$ and $d$ go to infinity. Therefore $[\mathbf{w}^\top \widehat{\mathbf{\Sigma}} \mathbf{w} - \widehat{U}(\gamma), \mathbf{w}^\top \widehat{\mathbf{\Sigma}} \mathbf{w} + \widehat{U}(\gamma)]$ is a valid level $(1-\gamma)100\%$ interval covering the true $\mathbf{w}^\top \mathbf{\Sigma} \mathbf{w}$.

## 2.2 Additional Data

This section considers the setting that there are available historical data for estimating $\mathbf{D}$. To adapt to the current market condition, we usually pick a short time series such that the asset returns are approximately stationary. However, it is likely that each univariate time series is stationary over a longer time scale than the multivariate time series, and hence we can incorporate extra information into calculation of the marginal standard deviations.

Inspired by this, we consider a setting where historical information is available. We do not assume the historical data to be multivariately stationary, but only marginally stationary. Formally speaking, let $\boldsymbol{R}_1, \ldots, \boldsymbol{R}_T$ be the observed stationary multivariate time series, and $\boldsymbol{H}_1, \ldots, \boldsymbol{H}_{T_h}$ be the available historical data with $\boldsymbol{H}_t = (H_{t1}, \ldots, H_{td})^\mathsf{T}$ and

$$T = O(T_h^{1-\delta}), \quad \text{where } \delta \text{ is an absolute constant.} \tag{2.9}$$

$\boldsymbol{H}_1, \ldots, \boldsymbol{H}_{T_h}$ could have overlap with $\boldsymbol{R}_1, \ldots, \boldsymbol{R}_T$. However, $\boldsymbol{H}_t$ is not necessarily identically distributed to either $\boldsymbol{H}_{t'}$ or $\boldsymbol{R}_1$ for any $t \neq t' \in \{1, \ldots, T_h\}$. Instead, we only assume that

$$H_{1j} \overset{\mathsf{d}}{=} H_{2j} \overset{\mathsf{d}}{=} \cdots \overset{\mathsf{d}}{=} H_{T_h j} \text{ and } \operatorname{Var}(H_{1j}) = \operatorname{Var}(R_{1j}), \quad \text{for } j \in \{1, \ldots, d\}.$$

We then estimate $\mathbf{w}^\mathsf{T}\boldsymbol{\Sigma}\mathbf{w}$ by separately estimating $\mathbf{D}$ and $\boldsymbol{\Sigma}^0$.

Formally, for estimating $\mathbf{D}$, we use the historical data $\boldsymbol{H}_1, \ldots, \boldsymbol{H}_{T_h}$ and derive

$$\widehat{\mathbf{D}}^h = (\widehat{\mathbf{D}}^h_{11}, \ldots, \widehat{\mathbf{D}}^h_{dd}), \quad \text{where } \widehat{\mathbf{D}}^h_{jj} := \widehat{\sigma}^h_{\mathrm{M},j} \frac{\widehat{\sigma}^h_1}{\widehat{\sigma}^h_{\mathrm{M},1}}, \tag{2.10}$$

and $\widehat{\sigma}^h_{\mathrm{M},j} = \widehat{\sigma}_{\mathrm{M}}(\{H_{tj}\}_{t=1}^T)$, for $j = 1, \ldots, d$, is the median absolute deviation estimator of $\{H_{tj}\}_{t=1}^T$, and $\widehat{\sigma}^h_1 = \big(\widehat{\operatorname{Var}}(\{H_{t1}\}_{t=1}^T)\big)^{1/2}$ is the Pearson sample standard deviation of $\{H_{t1}\}_{t=1}^T$. For estimating $\boldsymbol{\Sigma}^0$, we calculate the Kendall's tau matrix $\widehat{\mathbf{T}}$ based on $\{\boldsymbol{R}_1, \ldots, \boldsymbol{R}_T\}$.

**Remark 2.1.** In (2.10), to calculate $\widehat{\mathbf{D}}^h$, we employ the term $\widehat{\sigma}^h_1/\widehat{\sigma}^h_{\mathrm{M},1}$ to approximate the scaling factor between the median absolute deviation and the Pearson's standard deviation. This facilitates theoretical derivations. In practice, we can use, for example, the average version $\sum_{j=1}^d \widehat{\sigma}^h_j / \sum_{j=1}^d \widehat{\sigma}^h_{\mathrm{M},j}$ to estimate the scaling factor.

For estimating $\mathbf{w}^\mathsf{T}\boldsymbol{\Sigma}\mathbf{w}$, we replace $\widehat{\mathbf{D}}_{\text{est}}$ by $\widehat{\mathbf{D}}^h$, $\widehat{\boldsymbol{\Sigma}}^0_{\text{est}}$ by $\sin(\pi\widehat{\mathbf{T}}/2)$, and $\widehat{\boldsymbol{\Sigma}}_{\text{est}}$ by $\widehat{\boldsymbol{\Sigma}}^h$ in (2.6). For any given $1 - \gamma \in (0, 1)$, we calculate the Robust H-CLUB estimator $\widehat{U}^h(\gamma)$ as

$$\widehat{U}^h(\gamma) = \Phi^{-1}(1 - \gamma/2)\sqrt{\widehat{\sigma}^2_h/T}, \tag{2.11}$$

where $\widehat{\sigma}^2_h$ is calculated by employing the circular block bootstrap method introduced earlier. The corresponding confidence interval for the risk is

$$\big[\mathbf{w}^\mathsf{T}\widehat{\boldsymbol{\Sigma}}^h\mathbf{w} - \widehat{U}^h(\gamma), \mathbf{w}^\mathsf{T}\widehat{\boldsymbol{\Sigma}}^h\mathbf{w} + \widehat{U}^h(\gamma)\big]. \tag{2.12}$$

## 2.3 Unknown Marginal Volatilities

This section considers the setting that $\mathbf{D}$ is unknown with no additional data available. More precisely, we use a data splitting strategy for separately estimating $\mathbf{D}$ and $\boldsymbol{\Sigma}^0$. More precisely, we estimate $\mathbf{D}$ using the whole dataset:

$$\widehat{\mathbf{D}} = (\widehat{\mathbf{D}}_{11}, \ldots, \widehat{\mathbf{D}}_{dd}), \quad \text{with } \widehat{\mathbf{D}}_{jj} := \widehat{\sigma}_{\mathrm{M},j} \frac{\widehat{\sigma}_1}{\widehat{\sigma}_{\mathrm{M},1}}, \tag{2.13}$$

where $\widehat{\sigma}_{\mathrm{M},j} = \widehat{\sigma}_{\mathrm{M}}(\{R_{tj}\}_{t=1}^T)$ for $j = 1, \ldots, d$ and $\widehat{\sigma}_1 = \big(\widehat{\mathrm{Var}}(\{R_{t1}\}_{t=1}^T)\big)^{1/2}$ is the Pearson sample standard deviation of $\{R_{t1}\}_{t=1}^T$. For estimating $\mathbf{\Sigma}^0$, we extract a subsequence $\boldsymbol{R}_{T-T_s+1}, \ldots, \boldsymbol{R}_T$ from the time series $\boldsymbol{R}_1, \ldots, \boldsymbol{R}_T$, where $T_s \asymp T^{1-\delta}$ with $\delta$ a small enough absolute constant. Using this subsequence, we calculate the Kendall's matrix $\widehat{\mathbf{T}}^s$. Combining it with $\widehat{\mathbf{D}}$, we obtain a robust covariance matrix estimator

$$\widehat{\mathbf{\Sigma}}^s := \widehat{\mathbf{D}} \sin\Big(\frac{\pi}{2}\widehat{\mathbf{T}}^s\Big)\widehat{\mathbf{D}}.$$

We then estimate $\mathbf{w}^\mathsf{T}\mathbf{\Sigma}\mathbf{w}$ via replacing $\widehat{\mathbf{D}}_{\mathrm{est}}$, $\widehat{\mathbf{\Sigma}}^0_{\mathrm{est}}$, and $\widehat{\mathbf{\Sigma}}_{\mathrm{est}}$ by $\widehat{\mathbf{D}}$, $\sin(\frac{\pi}{2}\widehat{\mathbf{T}}^s)$, and $\widehat{\mathbf{\Sigma}}^s$ in (2.6). We then obtain a Robust H-CLUB estimator as

$$\widehat{U}^s(\gamma) = \Phi^{-1}(1-\gamma/2)\sqrt{\widehat{\sigma}_s^2/T_s}, \tag{2.14}$$

where $\widehat{\sigma}_s^2$ is calculated by employing the circular block bootstrap method. Accordingly, we construct the confidence interval of the risk as

$$\big[\mathbf{w}^\mathsf{T}\widehat{\mathbf{\Sigma}}^s\mathbf{w} - \widehat{U}^s(\gamma), \mathbf{w}^\mathsf{T}\widehat{\mathbf{\Sigma}}^s\mathbf{w} + \widehat{U}^s(\gamma)\big]. \tag{2.15}$$

**Remark 2.2.** In (2.13), for estimating the scaling factor, we can employ a similar average version as in Remark 2.1. We also note that the data splitting strategy is mainly proposed for theoretical analysis. In practice, we can set $\delta = 0$ and use the entire data set in calculating $\widehat{\mathbf{\Sigma}}^s$ and performing the block bootstrap.

# 3 Asymptotic Theory

In this section we prove that the confidence intervals of $\mathbf{w}^\mathsf{T}\mathbf{\Sigma}\mathbf{w}$ corresponding to three settings discussed in Section 2 have desired coverage probability. In other words, we prove that the Robust H-CLUB estimators proposed in (2.7), (2.11), and (2.14) are asymptotic $(1-\gamma)100\%$ confidence upper bound for the risk. It is clear that this problem reduces to calculating the limiting distributions of $\mathbf{w}^\mathsf{T}(\widehat{\mathbf{\Sigma}}_{\mathrm{est}} - \mathbf{\Sigma})\mathbf{w}$ for $\widehat{\mathbf{\Sigma}}_{\mathrm{est}} = \widehat{\mathbf{\Sigma}}, \widehat{\mathbf{\Sigma}}^h$, and $\widehat{\mathbf{\Sigma}}^s$. In the sequel, we adopt the triangular array setting as in Fan and Peng (2004) and Greenshtein and Ritov (2004) and allow the dimension $d$ to increase with the sample size $n$.

We introduce several mixing conditions for measuring degree of dependence. We start with an introduction of three mixing coefficients. For a $d$-dimensional stationary process $\{\boldsymbol{R}_t\}_{t\in\mathbb{Z}}$, let $\mathcal{F}_a^b$ be the $\sigma$-algebra generated by $\boldsymbol{R}_a, \ldots, \boldsymbol{R}_b$ for $a \le b$. We define the $\alpha$-, $\beta$-, and $\phi$-mixing coefficients as follows:

$$\alpha(n) := \sup_{B\in\mathcal{F}^0_{-\infty}, A\in\mathcal{F}^\infty_n} \big|\mathbb{P}(A\cap B) - \mathbb{P}(A)\mathbb{P}(B)\big|,$$

$$\beta(n) := \mathbb{E}\Big\{\sup_{A\in\mathcal{F}^\infty_n} \big|\mathbb{P}(A|\mathcal{F}^0_{-\infty}) - \mathbb{P}(A)\big|\Big\},$$

$$\phi(n) := \sup_{B\in\mathcal{F}^0_{-\infty}, A\in\mathcal{F}^\infty_n, \mathbb{P}(B)>0} \big|\mathbb{P}(A|B) - \mathbb{P}(A)\big|.$$

For an arbitrary positive integer $n$, we have $\alpha(n) \le \beta(n) \le \phi(n)$ (Yoshihara, 1976).

Suppose that $\{\boldsymbol{R}_1, \ldots, \boldsymbol{R}_T\}$ is a subsequence of the stationary process $\{\boldsymbol{R}_t\}_{t\in\mathbb{Z}}$. Let $F$ be the distribution function of $\boldsymbol{R}_1$. For $\boldsymbol{a} := \mathbf{D}\mathbf{w} = (a_1, \ldots, a_d)^\intercal$, let $g : \mathbb{R}^d \times \mathbb{R}^d \to \mathbb{R}$ be a kernel function

$$g(\boldsymbol{R}_t, \boldsymbol{R}_{t'}) := \frac{\pi}{2}\sum_{j\neq k} a_j a_k \cos(\frac{\pi}{2}\tau_{jk})\text{sign}(R_{tj} - R_{t'j})\text{sign}(R_{tk} - R_{t'k}). \tag{3.1}$$

We further define the following 3 quantities which will be useful in the later sections:

$$g_1(\boldsymbol{R}_1) := \int g(\boldsymbol{R}_1, \boldsymbol{R}_2)dF(\boldsymbol{R}_2), \tag{3.2}$$

$$\theta := \int g(\boldsymbol{R}_1, \boldsymbol{R}_2)dF(\boldsymbol{R}_1)dF(\boldsymbol{R}_2) = \boldsymbol{a}^\intercal\Big\{\cos(\frac{\pi}{2}\mathbf{T}) \circ \frac{\pi}{2}\mathbf{T}\Big\}\boldsymbol{a}, \tag{3.3}$$

$$\sigma^2 := 4\Big(\mathbb{E}g_1(\boldsymbol{R}_1)^2 - \theta^2 + 2\sum_{h=1}^{\infty}\Big\{\mathbb{E}g_1(\boldsymbol{R}_1)g_1(\boldsymbol{R}_{1+h})\Big\}\Big). \tag{3.4}$$

In the following, we assume that the elliptical time series model in Section 2 holds.

## 3.1 Theory for Known Volatilities

We make the following four assumptions which regulate the portfolio allocation vector $\mathbf{w}$ and the stationary process $\{\boldsymbol{R}_t\}_{t\in\mathbb{Z}}$.

**(A1)** There exist absolute constants $C_1$ and $C_2$ such that $\|\mathbf{w}\|_1 \le C_1$ and $\|\boldsymbol{\Sigma}\|_{\max} \le C_2$.

**(A2)** $\sigma$ is lower bounded by a positive absolute constant.

**(A3)** The process $\{\boldsymbol{R}_t\}_{t\in\mathbb{Z}}$ is $\phi$-mixing with $\phi(n) \le n^{-1-\epsilon}$ for some $\epsilon > 0$.

**(A4)** $\log d/(T^{1/2}) = o(1)$.

Assumption **(A1)** regulates the portfolio allocation vector $\mathbf{w}$ to prevent extreme positions. It is a common assumption made for stability of the portfolio (Jagannathan and Ma, 2003;

Fan et al., 2012, 2015). Assumption **(A2)** guarantees that the portfolio risk can not be diversified away. This is mild given that the returns are commonly assumed to follow a factor model (Chamberlain, 1983; Fan et al., 2015). Assumption **(A3)** is routinely used in analyzing time series to capture the serial dependence strength (Pan and Yao, 2008; Han and Liu, 2013b). Lastly, Assumption **(A4)** allows $d$ to grow nearly exponentially faster than $T$ and hence is mild.

In the setting of Section 2.1 and Assumptions **(A1)-(A4)**, we derive the limiting distribution of $\mathbf{w}^\mathsf{T}(\widehat{\mathbf{\Sigma}}-\mathbf{\Sigma})\mathbf{w}$. The following theorem shows that $\sqrt{T}\mathbf{w}^\mathsf{T}(\widehat{\mathbf{\Sigma}}-\mathbf{\Sigma})\mathbf{w}/\sigma$ is asymptotically normal.

**Theorem 3.1** (CLT, known volatilities)**.** Assuming that **(A0)** - **(A4)** hold and in the setting of Section 2.1, we have

$$\sqrt{T}\mathbf{w}^\mathsf{T}(\widehat{\mathbf{\Sigma}}-\mathbf{\Sigma})\mathbf{w}/\sigma \xrightarrow{\mathsf{d}} N(0,1),$$

as both $T$ and $d$ go to infinity.

The following theorem verifies that $\widehat{\sigma}^2$ calculated using the circular block bootstrap approach is a consistent estimator of $\sigma^2$. This result, combined with Theorem 3.1 and Slutsky's theorem, confirms that $\sqrt{T}\mathbf{w}^\mathsf{T}(\widehat{\mathbf{\Sigma}}-\mathbf{\Sigma})\mathbf{w}/\widehat{\sigma}$ converges weakly to the standard Gaussian. Accordingly, the confidence interval in (2.8) gives a reliable coverage probaility.

**Theorem 3.2** (bootstrap, known volatilities)**.** Under Assumptions **(A0)** - **(A4)**, we have

$$\widehat{\sigma}^2 = \sigma^2\big(1+o_P(1)\big),$$

and accordingly, for any given $\gamma \in (0,1)$, as $T,d\to\infty$, we have

$$\mathbb{P}\Big(\mathbf{w}^\mathsf{T}\mathbf{\Sigma}\mathbf{w} \in \big[\mathbf{w}^\mathsf{T}\widehat{\mathbf{\Sigma}}\mathbf{w} - \widehat{U}(\gamma), \mathbf{w}^\mathsf{T}\widehat{\mathbf{\Sigma}}\mathbf{w} + \widehat{U}(\gamma)\big]\Big) \to 1-\gamma.$$

The above two theorems only assume that the marginal second moments exist. Therefore, the Robust H-CLUB estimator naturally handles heavy-tailed data.

## 3.2 Theory with Additional Data

In this section we study the setting in Section 2.2. When $\mathbf{D}$ is unknown, we require additional assumptions. First, the following three assumptions require that $d$ does not grow too fast compared to $n$ and the given time series $\{\boldsymbol{X}_t\}_{t\in\mathbb{Z}}$ (either $\{\boldsymbol{R}_t\}_{t\in\mathbb{Z}}$ or $\{\boldsymbol{H}_t\}_{t\in\mathbb{Z}}$) is $\phi$-mixing with an exponentially decaying serial dependence.

- **(A5).** $\max\{\sqrt{\log d/T^\delta}, \log d/(T^{1/2})\} = o(1)$.

- **(A6).** The process $\{\boldsymbol{X}_t\}_{t\in\mathbb{Z}}$ is $\phi$-mixing with $\phi(n) \leq C_1 \exp(-C_2 n^r)$ for some absolute constants $C_1, C_2, r > 0$.
- **(A7).** Letting $a = \max(1, 1/r)$, we require that $\log d = o(T^{1/(2a+3)})$.

Recall that $\delta$ is defined in (2.9) for characterizing the length of historical data. Secondly, we require that the returns' $(4 + \epsilon_1)$-th moments exist for some absolute constant $\epsilon_1 > 0$, and the density functions are bounded away from zero around the median:

- **(A8).** For any $j \in \{1, \ldots, d\}$, $\mathbb{E}|X_{1j}|^{4+\epsilon_1} \leq C_0 < \infty$ for some constant $\epsilon_1, C_0 > 0$.
- **(A9).** Let $f_j$ and $\bar{f}_j$ be the density functions of $X_j$ and $|X_j - Q(X_j; 1/2)|$. For any $j \in \{1, \ldots, d\}$, we require $\inf_{|x-Q(f;1/2)|<\kappa} f(x) \geq \eta$ for some positive absolute constants $\kappa$ and $\eta$, and any $f \in \{f_j, \bar{f}_j\}$.

Under **(A0)** - **(A2)** and **(A5)** - **(A9)**, the next theorem shows that $\sqrt{T}\mathbf{w}^\mathsf{T}(\widehat{\mathbf{\Sigma}}^h - \mathbf{\Sigma})\mathbf{w}$ is asymptotically normal.

**Theorem 3.3** (CLT, unknown volatilities with additional data)**.** Assume that Assumptions **(A0)** - **(A2)** hold. In addition, assume that Assumptions **(A5)** - **(A7)** hold for both $\{\boldsymbol{R}_t\}_{t\in\mathbb{Z}}$ and the additional data $\{\boldsymbol{H}_t\}_{t\in\mathbb{Z}}$, and Assumptions **(A8)** - **(A9)** hold for $\{\boldsymbol{H}_t\}_{t\in\mathbb{Z}}$. Then in the setting of Section 2.2, we have

$$\sqrt{T}\mathbf{w}^\mathsf{T}(\widehat{\mathbf{\Sigma}}^h - \mathbf{\Sigma})\mathbf{w}/\sigma \xrightarrow{\mathsf{d}} N(0,1),$$

as both $T$ and $d$ go to infinity.

The next theorem shows that $\widehat{\sigma}_h^2$ is a consistent estimator of $\sigma^2$ and accordingly the confidence interval in (2.12) is valid.

**Theorem 3.4** (bootstrap, unknown volatilities with additional data)**.** Under the assumptions of Theorem 3.3, we have

$$\widehat{\sigma}_h^2 = \sigma^2\{1 + o_P(1)\},$$

and accordingly, for any given $\gamma \in (0,1)$, as $T, d \to \infty$, we have

$$\mathbb{P}\Big(\mathbf{w}^\mathsf{T}\mathbf{\Sigma}\mathbf{w} \in \big[\mathbf{w}^\mathsf{T}\widehat{\mathbf{\Sigma}}^h\mathbf{w} - \widehat{U}^h(\gamma), \mathbf{w}^\mathsf{T}\widehat{\mathbf{\Sigma}}^h\mathbf{w} + \widehat{U}^h(\gamma)\big]\Big) \to 1 - \gamma.$$

## 3.3 Theory with Unknown Marginal Volatilities

Lastly we study the setting in Section 2.3. Under this setting, we use a data splitting strategy and make inference only on a subsequence of length $T^{1-\delta}$. The next theorem justifies the use of such an approach.

**Theorem 3.5** (CLT, unknown marginal volatilities)**.** Assume that Assumptions **(A0)** - **(A2)** hold and Assumptions **(A5)** - **(A9)** hold for $\{\boldsymbol{R}_t\}_{t\in\mathbb{Z}}$. Then, under the setting of Section 2.3, we have

$$\sqrt{T_s}\mathbf{w}^\intercal(\widehat{\mathbf{\Sigma}}^s - \mathbf{\Sigma})\mathbf{w}/\sigma \xrightarrow{\mathsf{d}} N(0,1).$$

Furthermore, the bootstrap-based estimator $\widehat{\sigma}_s^2$ proves to be a consistent estimator of $\sigma^2$.

**Theorem 3.6** (bootstrap, unknown marginal volatilities)**.** Under the assumptions of Theorem 3.5, we have

$$\widehat{\sigma}_s^2 = \sigma^2\{1 + o_P(1)\},$$

and accordingly, for any given $\gamma \in (0,1)$, as $T, d \to \infty$, we have

$$\mathbb{P}\Big(\mathbf{w}^\intercal\mathbf{\Sigma}\mathbf{w} \in \big[\mathbf{w}^\intercal\widehat{\mathbf{\Sigma}}^s\mathbf{w} - \widehat{U}^s(\gamma), \mathbf{w}^\intercal\widehat{\mathbf{\Sigma}}^s\mathbf{w} + \widehat{U}^s(\gamma)\big]\Big) \to 1 - \gamma.$$

**Remark 3.7.** Compared to the method in Fan et al. (2015), the Robust H-CLUB estimator gains substantial robustness since it only assumes that the $(4+\epsilon_1)$-th moments exist for the marginal returns. In comparison, Fan et al. (2015) require a strong exponentially decaying rate in the tails (Check, for example, Assumption 3.4 therein). Such assumptions are often too restrictive and rarely satisfied in real applications. The Robust H-CLUB estimator attains the power for handling heavy-tailed data at the cost of a small $T^\delta$ efficiency. This is due to the data splitting strategy, which is an artifact of the proof. In practice, we find that the method introduced in Section 2.3 performs well.

The data splitting strategy allows the portfolio allocation vector to be random. More specifically, suppose that $\widehat{\mathbf{w}}$ is calculated based on the data $\boldsymbol{R}_1, \ldots, \boldsymbol{R}_T$. The next theorem shows that $\sqrt{T_s}\widehat{\mathbf{w}}^\intercal(\widehat{\mathbf{\Sigma}}^s - \mathbf{\Sigma})\widehat{\mathbf{w}}$ is asymptotically normal under assumptions outlined below.

**Corollary 3.1.** Under the assumptions in Theorem 3.5, let $\widehat{\mathbf{w}} = (\widehat{w}_1, \ldots, \widehat{w}_d)^\intercal$ be an estimator of $\mathbf{w} = (w_1, \ldots, w_d)^\intercal$ satisfying that

$$\mathbb{P}(|\widehat{w}_j/w_j - 1| > t) \le 2\exp\big(-CTt^2\big) \tag{3.5}$$

for some absolute constant $C$, any $j \in \{1, \ldots, d\}$, and any $t > 0$. We then have, as $T, d \to \infty$,

$$\sqrt{T_s}\widehat{\mathbf{w}}^\intercal(\widehat{\mathbf{\Sigma}}^s - \mathbf{\Sigma})\widehat{\mathbf{w}}/\sigma \xrightarrow{\mathsf{d}} N(0,1).$$

In this case, we can also employ a similar circular block bootstrap procedure for estimating the asymptotic variance of $\sqrt{T_s}\widehat{\mathbf{w}}^\intercal(\widehat{\mathbf{\Sigma}}^s - \mathbf{\Sigma})\widehat{\mathbf{w}}$.

# 4 Simulations on Synthetic Data

In this section we examine the finite-sample performance of the Robust H-CLUB estimators on synthetically generated data with heavy tails and noise contamination. We calculate several statistics of the estimators, following those used in Fan et al. (2015), to show the quality of the estimators. Our analysis shows that the Robust H-CLUB estimator performs well in all of the cases considered when compared to the full-confidence bound $\xi_T = \|\mathbf{w}\|_1^2 \|\widehat{\boldsymbol{\Sigma}}_{\text{est}} - \boldsymbol{\Sigma}\|_{\max}$. We observe that 95% confidence intervals by our proposed method are much tighter than the bound given by $\xi_T$. We also demonstrate that the H-CLUB calculated based on the robust estimators outperforms the H-CLUB based on the sample covariance matrix estimator $\mathbf{S}$ proposed in Fan et al. (2012) in the presence of heavy-tailed data. In particular, we show that the H-CLUB estimator does not achieve coverage proportions of 95% in the heavy-tailed setting, while the performance of the Robust H-CLUB estimator is consistently reliable. Lastly, we show that the Robust H-CLUB estimators also perform competitively when applied to the Gaussian data.

## 4.1 Calibration and Parameter Selection

To calibrate the parameters governing data generation in our model, we use the daily returns of the S&P 500's top 100 stocks ranked by market capitalization (as of June 29th, 2012), and the 3-month Treasury bill rates, sourced from the COMPUSTAT database (www.compustat.com) and the CSRP database (www.crsp.com), respectively. We consider the excess returns $\{\widetilde{\mathbf{y}}_t\}$ over the period from July 1, 2008 to June 29, 2012. We extract the following features:

1. $\{d_i^\dagger\}_{i=1}^{100}$ with $d_i^\dagger$ equal to the sample standard deviation of the $i$-th stock.
2. $\boldsymbol{\Sigma}^{0\dagger} = \{\Sigma_{ij}^{0\dagger}\}_{i,j=1}^{100}$, the sample correlation matrix of the observations $\widetilde{\mathbf{y}}_t$.

From these, we extract the mean and variance of $\{d_i^\dagger\}_{i=1}^{100}$, denoted respectively by $\mu_{\mathbf{d}^\dagger}$ and $\sigma_{\mathbf{d}^\dagger}^2$. We also compute the average and standard deviation of all pairwise correlations, denoted respectively by $\mu_{\boldsymbol{\Sigma}^{0\dagger}}$ and $\sigma_{\boldsymbol{\Sigma}^{0\dagger}}^2$. These parameters are used to generate correlation matrices and marginal variances later on.

We also have several tuning parameters to select. We choose $T_h = \lceil T^{1/(1-\delta_h)} \rceil$ with $\delta_h = 0.1$ as the parameter determining the quantity of historical data available to the estimator $\widehat{\boldsymbol{\Sigma}}^h$, $l = \lfloor T^{1-\epsilon_0} \rfloor$ with $\epsilon_0 = 0.8$ as the parameter controlling the block size in the block bootstrap, $N_{\text{bootstrap}} = 50$ as the number of bootstrapped datasets generated, and $T_s = \lfloor T^{1-\delta} \rfloor$ with $\delta = 0.01$ as the parameter controlling the data-splitting used in the estimator $\widehat{\boldsymbol{\Sigma}}^s$.

## 4.2 Simulation

For each given gross exposure constraint $c := \|\mathbf{w}\|_1$, we set $T = 300$ and allow $d$ to range from 50 to 500 in multiples of 50. For each value of $d$ we conduct 200 iterations of the same procedure: Generate a model, synthesize data from that model, and then calculate estimates based on the synthesized data. We collate the outputs across these 200 iterations to allow us to compare performance between different estimators.

The detailed procedure is described as follows:

1. Generate $\{d_i\}_{i=1}^d$ independently from the Gamma distribution with mean $\mu_{\mathbf{d}^\dagger}$ and variance $\sigma^2_{\mathbf{d}^\dagger}$. Define $\mathbf{D}$ as the diagonal matrix such that $\mathbf{D}_{ii} = d_i$.

2. Generate entries $\{\Sigma_{ij}\}_{i\neq j}$ of $\mathbf{\Sigma}^0$ independently from the Gaussian distribution with mean $\mu_{\mathbf{\Sigma}^{0\dagger}}$ and variance $\sigma^2_{\mathbf{\Sigma}^{0\dagger}}$. We threshold these off-diagonal elements to be no greater than 0.95 and set the diagonals of $\mathbf{\Sigma}^0$ to be 1. If the matrix is not positive definite, we use Higham's algorithm (see, e.g. Higham (2002)) to make it so, while keeping the diagonals fixed at 1.

3. Define the covariance matrix $\mathbf{\Sigma} = \mathbf{D}\mathbf{\Sigma}^0\mathbf{D}$.

4. Generate $\{\mathbf{R}_t\}_{t=1}^T$ independently from the multivariate $t$ distribution with 5 degrees of freedom and covariance matrix $\mathbf{\Sigma}$. Generate independent historical data $\{\mathbf{H}_t\}_{t=1}^{T_h}$ from the multivariate $t$ distribution with 5 degrees of freedom and covariance matrix $\mathbf{D}^2$.

5. Add noise contamination to the data by selecting a random 1% of the elements in $\{\mathbf{R}_t\}_{t=1}^T$ and multiplying each one by a random variable drawn independently from a Unif$(1, 15)$ distribution. Do the same to 1% of the elements in $\{\mathbf{H}_t\}_{t=1}^{T_h}$. This step can be regarded as the news arrivals on the firms that cause their returns to jump.

6. Calculate the covariance estimates given by the sample covariance matrix $\mathbf{S}$ and the robust estimators $\widehat{\mathbf{\Sigma}}$, $\widehat{\mathbf{\Sigma}}^h$, and $\widehat{\mathbf{\Sigma}}^s$, using the tuning parameters given in Section 4.1[4].

7. Generate 500 portfolio allocation vectors $\mathbf{w}$ according to the method outlined in Fan et al. (2015), which is approximately uniformly distributed on the manifold $\{\mathbf{w} : \|\mathbf{w}\|_1 = c, \mathbf{w}^T\mathbf{1} = 1\}$.

[4] We find the following minor alteration to improve performance in practice: For the H-CLUB based on $\widehat{\mathbf{\Sigma}}^h$, we take block-bootstrapped samples of both $\{\mathbf{H}_t\}_{t=1}^{T_h}$ and $\{\mathbf{R}_t\}_{t=1}^T$ in estimating the variance of $\mathbf{w}^\top\widehat{\mathbf{\Sigma}}^h\mathbf{w}$. For this we use the block size parameter $l_h = \lfloor T_h^{1-\epsilon_0}\rfloor$, entirely analogously to the block bootstrapping performed on $\{\mathbf{R}_t\}_{t=1}^T$ with $l = \lfloor T^{1-\epsilon_0}\rfloor$. We use this modification throughout Sections 4 and 5.

8. For each portfolio allocation, calculate the H-CLUB estimates corresponding to the estimators listed in Step 6. As proof-of-concept, we also calculate the estimator with $\widehat{\mathbf{\Sigma}}^s_{T_s=T}$, which is the estimator $\widehat{\mathbf{\Sigma}}^s$ with $T_s = T$ (i.e., no data-splitting performed).

9. Over the 500 portfolios, compute the averages of the true risk $R(\mathbf{w}) := \sqrt{\mathbf{w}^\mathsf{T}\mathbf{\Sigma}\mathbf{w}}$, as well as $\Delta := |\mathbf{w}^\mathsf{T}(\widehat{\mathbf{\Sigma}}_{\text{est}} - \mathbf{\Sigma})\mathbf{w}|, \xi_T := \|\mathbf{w}\|_1^2 \|\widehat{\mathbf{\Sigma}}_{\text{est}} - \mathbf{\Sigma}\|_{\max}$, and $\widehat{U}(0.05) = 2\sqrt{\widehat{\sigma}^2/T}$ for each of the estimators $\widehat{\mathbf{\Sigma}}_{\text{est}}$ considered.

We plot the averages of $\Delta$, $\xi_T$, and $\widehat{U}(0.05)$ against $d$ for every estimator considered and for $c = 1$, $c = 1.6$, and $c = 2$ to observe the effects of gross exposure on risk assessment.

Next, for $d = 200$ and $d = 500$, we calculate the following quantities over the 100,000 portfolios (500 portfolios over 200 synthetic datasets) : The coverage proportion, defined as the fraction of the sample in which the 95% confidence interval contains the true risk $R(\mathbf{w}) = (\mathbf{w}^\mathsf{T}\mathbf{\Sigma}\mathbf{w})^{1/2}$, the ratio of bounds defined as

$$\text{RE}_1 := \frac{\xi_T}{2\sqrt{\widehat{\sigma}^2/T}},$$

and the relative error defined as

$$\text{RE}_2 := \frac{\sqrt{\widehat{\sigma}^2/T}}{2\mathbf{w}^T\mathbf{\Sigma}\mathbf{w}}.$$

Again, we compute these for $c = 1$, 1.6, and 2. The measure $\text{RE}_1$ compares the upper bound with the half width of the 95% confidence interval, whereas $\text{RE}_2$ is the half width of 95% confidence interval for the portfolio risk $\{\mathbf{w}^T\mathbf{\Sigma}\mathbf{w}\}^{1/2}$ divided by the portfolio risk itself. The former depicts how inefficiency the confidence upper bound is and the latter measures how informative the constructed confidence interval is.

Lastly, we repeat the previous calculations of coverage proportions, $\text{RE}_1$ and $\text{RE}_2$ in a setting where the data are generated from a Gaussian distribution without any noise contamination. This means we alter Step 4 of the procedure above (but substitute Gaussian distribution for $t$ distribution) and remove Step 5. This allows us to examine the degree of efficiency loss for robustness when data are normal. In this setting, we also calculate the ratio $\widehat{U}(0.05)/\Delta$ as a measure of how tight the H-CLUB is relative to the theoretical minimum bound.

## 4.3 Results

In Figures 1 and 2, we plot the average risk estimation errors along with the estimated error bounds with gross exposure $c = 1$, 1.6, and 2, using estimators $\widehat{\mathbf{\Sigma}}_{\text{est}} = \widehat{\mathbf{\Sigma}}, \widehat{\mathbf{\Sigma}}^h, \widehat{\mathbf{\Sigma}}^s$, and $\widehat{\mathbf{\Sigma}}^s_{T_s=T}$. Note that $c = 1.6$ results in an average 130% long positions and 30% short positions, which is commonly used in practice. We also use the sample covariance matrix estimator $\mathbf{S}$, for which an H-CLUB estimator was derived in Fan et al. (2015), which is not robustified.

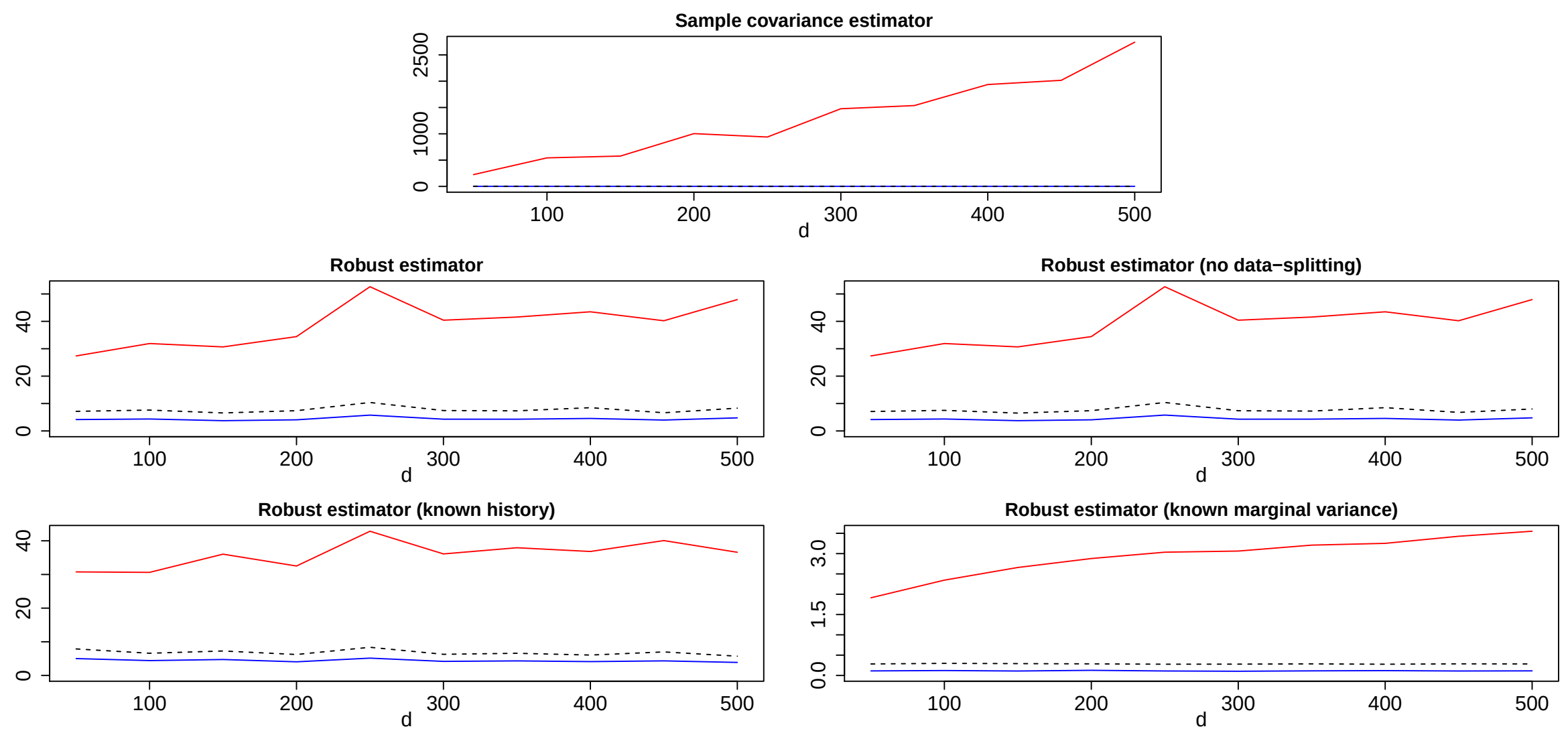

Figure 1: Averages of $\Delta = |\mathbf{w}^\intercal(\widehat{\mathbf{\Sigma}}_{\text{est}} - \mathbf{\Sigma})\mathbf{w}|$ (blue curve), $\widehat{U}(0.05) = 2\sqrt{\widehat{\text{Var}}(\mathbf{w}^\intercal\widehat{\mathbf{\Sigma}}_{\text{est}}\mathbf{w})}$ (dashed curve), and $\xi_T = \|\mathbf{w}\|_1^2\|\widehat{\mathbf{\Sigma}}_{\text{est}} - \mathbf{\Sigma}\|_{\max}$ (red curve) for $c = 1.0$. Horizontal axis shows dimension of problem, i.e., portfolio size. Vertical axis shows the calculated averaged values.

From these plots, we see that

- The dashed curve lies above the solid blue line throughout, an indication of the validity of the 95% bound given by $\widehat{U}(0.05)$. It is interesting to note that this still holds for the sample covariance matrix estimator $\mathbf{S}$, but this is in the average sense. As we will see in Table 1, however, $\mathbf{S}$ fails to attain 95% coverage.

- The crude bound $\xi_T$ is much larger than either the true error $\Delta$ or the 95% confidence bound $\widehat{U}(0.05)$. This discrepancy increases with $d$, but also with $c$ as we can see by comparing Figure 1 with Figure 2. This is quantified in Table 2.

- For large $d$ the crude bound on the sample covariance matrix estimator is almost 100 times larger than on any of the robust estimators. This suggests inaccurate estimation of the sample covariance in the presence of heavy tails and contamination.

Table 1 illustrates the coverage of each estimator, defined as the proportion of samples in which the 95% confidence interval captures the true variance $\mathbf{w}^\intercal\mathbf{\Sigma}\mathbf{w}$. It can be seen that all the robust estimators have coverage proportions of approximately 95%. However, the sample covariance matrix estimator $\mathbf{S}$ has substantially lower coverage. It is not sufficiently robust to give a valid bound under the current setting.

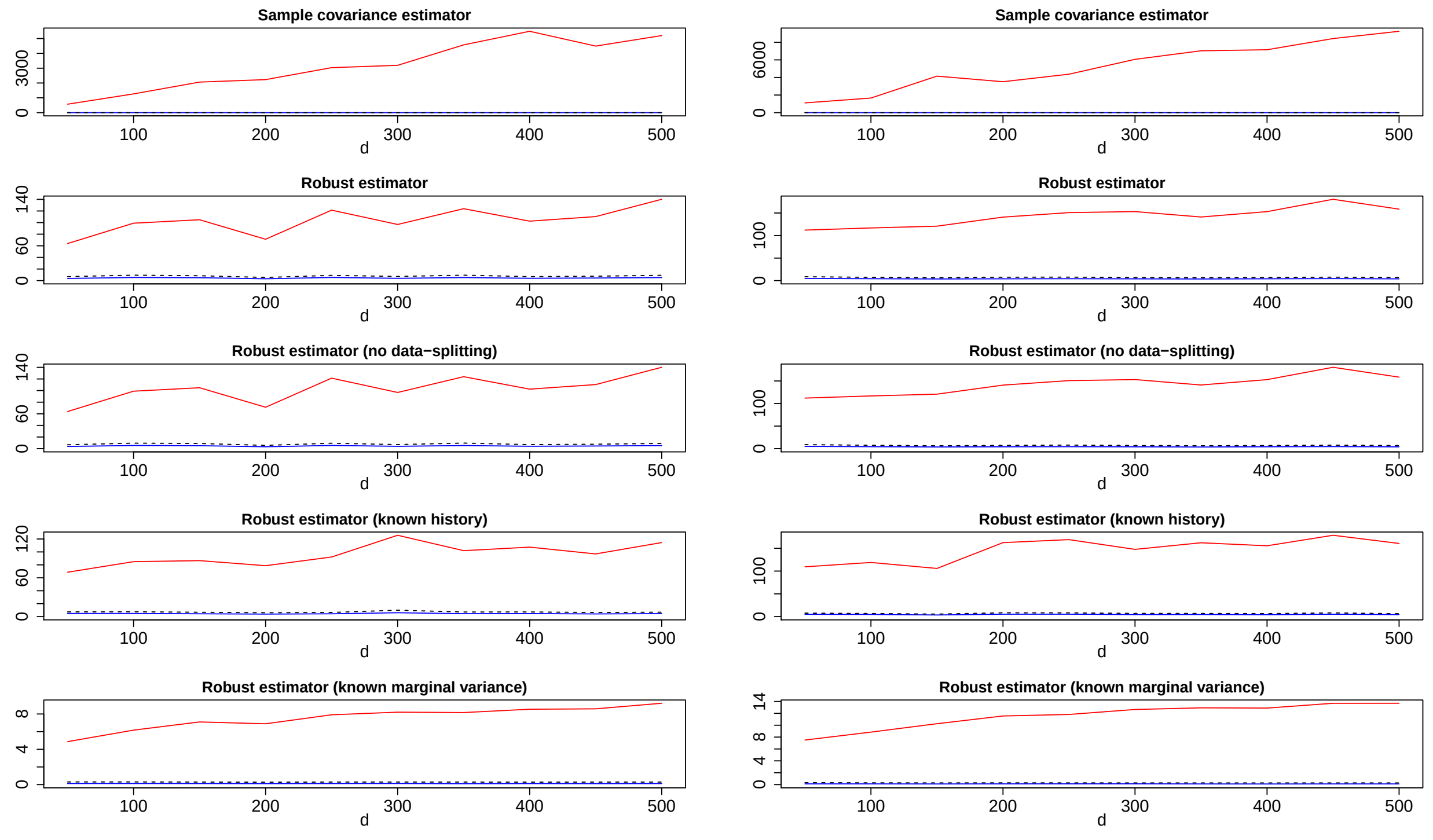

(a) $c = 1.6$ (b) $c = 2$

Figure 2: Averages of $\Delta = |\mathbf{w}^\top(\widehat{\mathbf{\Sigma}}_{\text{est}} - \mathbf{\Sigma})\mathbf{w}|$ (blue curve), $\widehat{U}(0.05) = 2\sqrt{\widehat{\text{Var}}(\mathbf{w}^\top \widehat{\mathbf{\Sigma}}_{\text{est}} \mathbf{w})}$ (dashed curve) and $\xi_T = \|\mathbf{w}\|_1^2 \|\widehat{\mathbf{\Sigma}}_{\text{est}} - \mathbf{\Sigma}\|_{\max}$ (red curve) for $c = 1.6$ and $c = 2$. Horizontal axis shows dimension of problem, i.e., portfolio size. Vertical axis shows the calculated averaged values.

We make further comparisons between the robust estimators we have proposed. Table 2 illustrates averages and standard deviations of the ratio $\text{RE}_1 = \xi_T/\widehat{U}(0.05)$: the ratio between the full confidence bound and the H-CLUB. These serve to quantify some of our observations made on Figures 1 and 2 — in particular, that the ratio $\xi_T/\widehat{U}(0.05)$ increases strongly with $c$ and weakly with $d$.

We observe that:

- The value of $\text{RE}_1$ is considerably bigger than 1, reflecting the fact that the confidence interval given by the Robust H-CLUB is much tighter than that given by the crude bound. In almost all cases the value of $\text{RE}_1$ reflects a difference of scale of an order of magnitude between the H-CLUB interval and the crude interval using $\xi_T$.

Table 1: Empirical coverage proportion for 95% confidence intervals in settings of data drawn from $t_5$ distribution with 1% noise contamination. Taken over 200 samples with $T = 300$.

| | $d = 200$ | | | $d = 500$ | | |
|---|---|---|---|---|---|---|
| | $c = 1.0$ | $c = 1.6$ | $c = 2.0$ | $c = 1.0$ | $c = 1.6$ | $c = 2.0$ |
| Coverage $\mathbf{S}$ | 81.88% | 72.29% | 69.31% | 83.30% | 82.24% | 80.12% |
| Coverage $\widehat{\mathbf{\Sigma}}^s$ | 97.59% | 95.26% | 97.64% | 99.00% | 97.09% | 95.52% |
| Coverage $\widehat{\mathbf{\Sigma}}^s_{T_s=T}$ | 96.38% | 95.70% | 97.49% | 98.18% | 97.03% | 95.03% |
| Coverage $\widehat{\mathbf{\Sigma}}^h$ | 93.87% | 93.19% | 95.23% | 93.01% | 92.84% | 94.67% |
| Coverage $\widehat{\mathbf{\Sigma}}$ | 94.21% | 95.54% | 96.40% | 95.16% | 93.41% | 93.67% |

- The ratio $\mathrm{RE}_1$ increases with our ability to accurately estimate the marginal standard deviations. Note that $\mathrm{RE}_1(\widehat{\mathbf{\Sigma}}) > \mathrm{RE}_1(\widehat{\mathbf{\Sigma}}^h) > \mathrm{RE}_1(\widehat{\mathbf{\Sigma}}^s_{T_s=T}) > \mathrm{RE}_1(\widehat{\mathbf{\Sigma}}^s)$, which corresponds to an ordering based on the amount of information used to estimate the marginal standard deviations.

- The value of $\mathrm{RE}_1$ increases strongly with $c$ and weakly with $d$. This suggests that the accuracy benefits of using the H-CLUB over the crude bound are particularly substantial for larger portfolios and those with higher gross exposure.

Table 3 summarizes the relative error ($\mathrm{RE}_2$), which shows how informative our confidence intervals for the true portfolio risks are. Similar to Table 2, we show the mean and standard deviation of $\mathrm{RE}_2$ calculated over 200 simulations with 500 randomly generated portfolios per simulation (i.e. 100,000 portfolios total).

Here we see a similar pattern as before. Values are generally better (smaller, here) when more information is available in our estimation of the marginal standard deviations. This statement comes from the observation that $\mathrm{RE}_2(\widehat{\mathbf{\Sigma}}) \ll \mathrm{RE}_2(\widehat{\mathbf{\Sigma}}^h) < \mathrm{RE}_2(\widehat{\mathbf{\Sigma}}^s_{T_s=T}) < \mathrm{RE}_2(\widehat{\mathbf{\Sigma}}^s)$. We also observe that here the value of $\mathrm{RE}_2$ does not appear to vary much with either $c$ or $d$. It is also substantially larger than the values seen in, e.g., Fan et al. (2015), presumably due to the heavier tails and presence of noise in the data here which is not seen in those settings. This difference can be immediately observed by comparing with Table 4. From the last row of Table 3, the uninformative construction of the confidence interval is mainly due to the inaccurate estimation of the marginal variances in presence of large random noises

Table 2: Averages and standard deviations (in parentheses) of $\mathrm{RE}_1 := \xi_T/(2\sqrt{\widehat{\sigma}^2/T})$ over 200 samples.

| | $d = 200$ | | | $d = 500$ | | |
|---|---|---|---|---|---|---|
| | $c = 1.0$ | $c = 1.6$ | $c = 2.0$ | $c = 1.0$ | $c = 1.6$ | $c = 2.0$ |
| $\mathrm{RE}_1$ | 5.57 | 14.73 | 21.62 | 6.63 | 17.55 | 27.50 |
| $\widehat{\mathbf{\Sigma}}^s$ | (1.94) | (5.51) | (7.68) | (2.18) | (6.13) | (9.95) |
| $\mathrm{RE}_1$ | 5.64 | 14.54 | 21.90 | 6.70 | 17.47 | 27.57 |
| $\widehat{\mathbf{\Sigma}}^s_{T_s=T}$ | (1.85) | (5.64) | (8.50) | (2.32) | (6.61) | (9.39) |
| $\mathrm{RE}_1$ | 5.87 | 14.65 | 22.44 | 6.93 | 18.54 | 27.22 |
| $\widehat{\mathbf{\Sigma}}^h$ | (2.11) | (5.24) | (8.55) | (2.25) | (6.56) | (9.55) |
| $\mathrm{RE}_1$ | 9.88 | 25.43 | 38.85 | 12.29 | 32.19 | 48.62 |
| $\widehat{\mathbf{\Sigma}}$ | (2.80) | (7.31) | (10.89) | (3.13) | (9.10) | (12.91) |

Table 3: Averages and standard deviations (in parentheses) of $\mathrm{RE}_2 = \sqrt{\widehat{\sigma}^2/T}/2\mathbf{w}^T\mathbf{\Sigma}\mathbf{w}$ over 200 samples.

| | $d = 200$ | | | $d = 500$ | | |
|---|---|---|---|---|---|---|
| | $c = 1.0$ | $c = 1.6$ | $c = 2.0$ | $c = 1.0$ | $c = 1.6$ | $c = 2.0$ |
| $\mathrm{RE}_2$ | 0.513 | 0.627 | 0.478 | 0.521 | 0.549 | 0.480 |
| $\widehat{\mathbf{\Sigma}}^s$ | (0.609) | (0.880) | (0.534) | (0.586) | (0.606) | (0.540) |
| $\mathrm{RE}_2$ | 0.500 | 0.644 | 0.483 | 0.517 | 0.559 | 0.471 |
| $\widehat{\mathbf{\Sigma}}^s_{T_s=T}$ | (0.594) | (0.906) | (0.554) | (0.595) | (0.626) | (0.531) |
| $\mathrm{RE}_2$ | 0.462 | 0.571 | 0.575 | 0.492 | 0.471 | 0.494 |
| $\widehat{\mathbf{\Sigma}}^h$ | (0.485) | (0.837) | (0.691) | (0.604) | (0.555) | (0.573) |
| $\mathrm{RE}_2$ | 0.022 | 0.021 | 0.021 | 0.021 | 0.021 | 0.021 |
| $\widehat{\mathbf{\Sigma}}$ | (0.002) | (0.002) | (0.002) | (0.002) | (0.002) | (0.002) |

and heavy tails.

For our last set of results on synthetic data, we show in Table 4 that the robust estimators are still competitive with the sample covariance based estimator when the data are drawn from a Gaussian distribution without noise contamination. In this table we present coverage proportions, means of $\mathrm{RE}_1$ and $\mathrm{RE}_2$, as well as the mean of the ratio between the 95% H-CLUB and the value it is upper bounding, with this ratio given by $\widehat{U}(0.05)/\Delta$. These are calculated over 200 randomly generated models.

Table 4: Coverage proportion and means of $\mathrm{RE}_1$, $\mathrm{RE}_2$ and $\widehat{U}(0.05)/\Delta$ for 200 samples when returns are drawn from Gaussian distributions without noise contamination, using $d = 500$.

| | Coverage | | | $\mathrm{RE}_1$ | | | $\mathrm{RE}_2$ | | | $\widehat{U}(0.05)/\Delta$ | | |
|---|---|---|---|---|---|---|---|---|---|---|---|---|
| $c$ | 1.0 | 1.6 | 2.0 | 1.0 | 1.6 | 2.0 | 1.0 | 1.6 | 2.0 | 1.0 | 1.6 | 2.0 |
| $\mathbf{S}$ | .948 | .944 | .927 | 8.10 | 21.22 | 33.17 | 4.01% | 3.97% | 4.01% | 5.67 | 6.29 | 7.02 |
| $\widehat{\mathbf{\Sigma}}^s$ | .965 | .954 | .950 | 8.57 | 22.24 | 34.19 | 7.19% | 7.13% | 7.14% | 7.88 | 5.86 | 7.09 |
| $\widehat{\mathbf{\Sigma}}^s_{T_s=T}$ | .960 | .951 | .950 | 8.58 | 22.46 | 33.98 | 7.17% | 7.06% | 7.20% | 7.42 | 5.93 | 8.88 |
| $\widehat{\mathbf{\Sigma}}^h$ | .960 | .953 | .964 | 9.26 | 23.99 | 37.28 | 4.92% | 5.01% | 5.09% | 6.97 | 7.14 | 5.47 |
| $\widehat{\mathbf{\Sigma}}$ | .957 | .949 | .923 | 11.65 | 30.65 | 48.75 | 2.01% | 2.00% | 2.00% | 7.05 | 6.76 | 5.82 |

# 5 An Empirical Study

In this section we examine the behaviour of the Robust H-CLUB estimators when applied to real-world data. We use the daily excess returns of 100 industrial portfolios formed on size and book-to-market ratio, as available on the website of Kenneth French. We use the subset of data spanning from July 1, 2008 to June 29, 2012. For each 21 day period (nominal month), we use the preceding 21 days' data to estimate the covariance matrix via the Robust H-CLUB estimator with data-splitting ($\widehat{\mathbf{\Sigma}}^s$), the Robust H-CLUB with no data splitting ($\widehat{\mathbf{\Sigma}}^s_{T_s=T}$), and the Robust H-CLUB estimator with known history ($\widehat{\mathbf{\Sigma}}^h$). For the matrix of additional observations used in the latter estimator, we use the preceding 1.5 months (31 days) of returns data. Note that for all robust estimators in this section we use the tuning parameter $l = \lfloor T^{0.5} \rfloor$ (i.e. $\epsilon_0 = 0.5$) for the block size in the bootstrapping procedure. All other parameters are as in the previous section. Finally, we also estimate the covariance via the sample covariance matrix estimator $\mathbf{S}$ for comparison.

We track the performance of the H-CLUB estimators on three portfolios: one portfolio with equal weighting ($\widehat{\mathbf{w}} = (1/100, \ldots, 1/100)$), and two portfolios of minimum variance with gross exposure $c = 1$ and $c = 1.6$, as given by

$$\widehat{\mathbf{w}} = \underset{\mathbf{w}^\mathsf{T}\mathbf{1}=1, \|\mathbf{w}\|_1 = c}{\arg\min} \ \mathbf{w}^\mathsf{T}\widehat{\mathbf{\Sigma}}_{\text{est}}\mathbf{w}.$$

Note that on occasion the estimated covariance matrix is not positive definite, leading to problems in solving for the portfolio of minimum variance. In these cases, we coerce the estimated covariance matrix to be positive definite using Higham's algorithm before calculating the minimum variance portfolio.

The portfolios of minimum variance are calculated at the start of each nominal month. The actual risk during the holding month for each $\widehat{\mathbf{w}}$ as defined above is then

$$R(\widehat{\mathbf{w}}) = (\widehat{\mathbf{w}}^\mathsf{T}\mathbf{\Sigma}\widehat{\mathbf{w}})^{1/2} \quad \text{and} \qquad \mathbf{\Sigma} = \frac{1}{21}\sum_{t=1}^{21}\mathbf{y}_t\mathbf{y}_t^\mathsf{T},$$

where $\{\mathbf{y}_t\}_{t=1}^T$ are the centralized daily returns over the holding month. This is calculated for each month in the four year period of study.

For each estimator and portfolio strategy, we consider five quantities. These quantities are summarized via their mean (calculated over the whole study period) in Table 5. We compare the first two columns of Table 5 and provide several observations.

- The values of $\Delta$ are comparable among the four estimators considered. This suggests that all estimators are similar in their estimations of the covariance matrix $\mathbf{\Sigma}$, and that differences between them lie in their ability to accurately conduct inference on $\widehat{\mathbf{\Sigma}}_{\text{est}}$ (i.e. construct a valid H-CLUB).

- The (non-robustified) sample covariance matrix estimator $\mathbf{S}$ fails to give a valid upper bound, as $\widehat{U}(0.05)$ is less than $\Delta$ throughout.

- For the robust estimators, $\widehat{U}(0.05)$ is greater than $\Delta$ for all cases except one. This is broadly consistent with the expectation that the value of $\widehat{U}(0.05)$ for the robust estimators is a 95% upper bound of the estimation error for portfolio variance. We note that for the single discrepancy ($\widehat{\mathbf{\Sigma}}^h$, on the minimum variance portfolio with $\|\mathbf{w}\|_1 = 1.6$), the value of $\widehat{U}(0.05)$ still only falls below $\Delta$ by a small margin.

Lastly, the estimated risk error $\widehat{U}(0.05)/\sqrt{4\mathbf{w}^\mathsf{T}\widehat{\mathbf{\Sigma}}_{\text{est}}\mathbf{w}}$ is an H-CLUB estimate for the true risk error $|(\mathbf{w}^\mathsf{T}\mathbf{\Sigma}\mathbf{w})^{1/2} - (\mathbf{w}^\mathsf{T}\widehat{\mathbf{\Sigma}}_{\text{est}}\mathbf{w})^{1/2}|$ (we can see this simply by applying the delta method to the results of, e.g. Theorem 3.6). The last two columns of Table 5 show that the robust estimators hold true to this, with the estimated risk error uniformly bounding the true risk error in all cases. However, the non-robustified sample covariance estimator does not yield a

Table 5: Annualized true and estimated risk errors calculated on the 100 Fama-French portfolios

| Strategy | Average of $\Delta(\times 10^{-4})$ | Average of $\widehat{U}(0.05)(\times 10^{-4})$ | Average of True Risk | True Risk Error | Estimated Risk Error |
|---|---|---|---|---|---|
| | | $\mathbf{S}$ *(Sample Covariance Matrix Estimator)* | | | |
| Equal weighted | 2.310 | 1.939 | 27.36% | 8.32% | 6.62% |
| Min. variance ($c = 1$) | 1.289 | 0.743 | 19.52% | 6.97% | 4.19% |
| Min. variance ($c = 1.6$) | 0.760 | 0.312 | 15.25% | 6.38% | 2.66% |
| | | $\widehat{\mathbf{\Sigma}}^s$ *(Robust Estimator)* | | | |
| Equal weighted | 2.165 | 4.790 | 27.36% | 8.35% | 18.67% |
| Min. variance ($c = 1$) | 1.470 | 2.696 | 21.06% | 8.41% | 17.67% |
| Min. variance ($c = 1.6$) | 1.576 | 2.249 | 18.30% | 13.05% | 46.32% |
| | | $\widehat{\mathbf{\Sigma}}^s_{T_s=T}$ *(Robust Estimator — no data-splitting)* | | | |
| Equal weighted | 2.154 | 5.121 | 27.36% | 8.32% | 18.94% |
| Min. variance ($c = 1$) | 1.459 | 2.826 | 21.02% | 8.34% | 20.41% |
| Min. variance ($c = 1.6$) | 1.562 | 2.218 | 18.22% | 12.86% | 37.81% |
| | | $\widehat{\mathbf{\Sigma}}^h$ *(Robust Estimator — known history)* | | | |
| Equal weighted | 2.100 | 3.325 | 27.36% | 7.69% | 12.85% |
| Min. variance ($c = 1$) | 1.390 | 1.885 | 20.79% | 7.63% | 12.25% |
| Min. variance ($c = 1.6$) | 1.358 | 1.200 | 17.52% | 10.99% | 17.40% |

*Note:* $\Delta = |\mathbf{w}^\top(\widehat{\mathbf{\Sigma}}_{\text{est}} - \mathbf{\Sigma})\mathbf{w}|$, $\widehat{U}(0.05) = 2 \times (\widehat{\text{Var}}(\mathbf{w}^\top \widehat{\mathbf{\Sigma}}_{\text{est}}\mathbf{w}))^{1/2}$. *True Risk is* $\sqrt{252} \times R(\mathbf{w})$*. True Risk Error is* $\sqrt{252} \times |(\mathbf{w}^\top \widehat{\mathbf{\Sigma}}_{\text{est}}\mathbf{w})^{1/2} - (\mathbf{w}^\top \mathbf{\Sigma}\mathbf{w})^{1/2}|$*, and Estimated Risk Error is* $\sqrt{252} \times \widehat{U}(0.05)/\sqrt{4\mathbf{w}^\top \widehat{\mathbf{\Sigma}}_{\text{est}}\mathbf{w}}$*. The factor of* $\sqrt{252}$ *is present to convert the risks to annualized values.*

good upper bound, with the estimated risk error uniformly falling below the true risk error. This is again an evidence for the strength of the proposed robust estimators in the presence of heavy-tailed or noisy data.

# 6 Conclusion and Discussion

This paper considers the problem of assessing the risks of large portfolios in a robust manner. We consider three different settings depending on whether $\mathbf{D}$ is known or not, and propose three corresponding Robust H-CLUB approaches based on robust rank and quantile statistics. For the first time in the literature, we provide an inferential theory of these robust risk estimators. Compared to Fan et al. (2015), the proposed approaches do not require strong moment assumptions on the data. Both theoretical and empirical results verify that the Robust H-CLUB approaches are more appropriate for studying heavy-tailed asset returns.

In the present paper, we do not impose any structural assumption on the covariance matrix, such as the low rank plus sparse structure induced by the factor model. Fan et al. (2015) propose methods based on factor-based covariance matrix estimators proposed in Fan et al. (2008) and Fan et al. (2013). A natural extension to Fan et al. (2013) is to use $\widehat{\mathbf{\Sigma}}$ (or $\widehat{\mathbf{\Sigma}}^h, \widehat{\mathbf{\Sigma}}^s$), instead of the sample covariance $\mathbf{S}$, as the pilot estimator and plug it into the POET algorithm (Fan et al., 2013). This constructs another robust risk estimator. We plan to investigate the theoretical properties of such robust risk estimators and their limiting distributions in the future.

The results in this paper also raise a number of interesting questions for future research. One example is on deriving the limiting distributions of functionals of $\widehat{\mathbf{\Sigma}}$ other than $\mathbf{w}^\mathsf{T}\widehat{\mathbf{\Sigma}}\mathbf{w}$. For example, Han and Liu (2014a) study the limiting distribution of $\|\widehat{\mathbf{\Sigma}}\|_{\max}$ as $T, d \to \infty$ in the setting that the observations are mutually independent. It is interesting to investigate such asymptotic theory for a multivariate time series.

# 7 Proofs

In this section we provide the proofs of results in Section 3. In the sequel, using Assumption **(A1)**, we assume that $\|\mathbf{w}\| = 1$ and $\|\mathbf{\Sigma}\|_{\max} \le 1$ without loss of generality.

## 7.1 Supporting Lemmas

**Lemma 7.1** (Kontorovich et al. (2008) and Mohri and Rostamizadeh (2010))**.** Let $f : \Omega^T \to \mathbb{R}$ be a measurable function that is $c$-Lipschitz with regard to the Hamming metric for some $c > 0$:

$$\sup_{x_1,\ldots,x_t,x_t'} \Big| f(x_1, \ldots, x_t, \ldots, x_T) - f(x_1, \ldots, x_t', \ldots, x_T) \Big| \le c,$$

and $X_1, \ldots, X_T$ be a sequence of stationary $\phi$-mixing random variables. Then, for any $\epsilon > 0$, the following inequality holds:

$$\mathbb{P}\Big\{ |f(X_1, \ldots, X_T) - \mathbb{E} f(X_1, \ldots, X_T)| \ge \epsilon \Big\} \le 2 \exp\Big[ - \frac{2\epsilon^2}{Tc^2\{1 + 2\sum_{k=1}^T \phi(k)\}} \Big].$$

**Lemma 7.2** (Yoshihara (1976))**.** Let $\{\boldsymbol{X}_t\}_{t\in\mathbb{Z}}$ be a stationary process with the distribution function $F$. For $T \ge m$, we define

$$U_T(g) = \binom{T}{m}^{-1} \sum_{i_1 < \cdots < i_m} g(\boldsymbol{X}_{i_1}, \ldots, \boldsymbol{X}_{i_m})$$

be a $U$-statistic with order $m$ and kernel function $g$. Let the function $g_i(\cdot)$ be defined as

$$g_i(\boldsymbol{X}_1, \ldots, \boldsymbol{X}_i) = \int g(\boldsymbol{X}_1, \ldots, \boldsymbol{X}_m) dF(\boldsymbol{X}_{i+1}) \ldots dF(\boldsymbol{X}_m),$$

for $1 \le i \le m$, and let parameters $\theta$ and $\sigma^2$ be defined as

$$\theta = \int g(\boldsymbol{X}_1, \ldots, \boldsymbol{X}_m) dF(\boldsymbol{X}_1) \ldots dF(\boldsymbol{X}_m),$$

$$\sigma^2 = 4\Big(\mathbb{E} g_1(\boldsymbol{X}_1)^2 - \theta^2 + 2\sum_{h=1}^{\infty}\Big(\mathbb{E} g_1(\boldsymbol{X}_1) g_1(\boldsymbol{X}_{1+h}) - \theta^2\Big)\Big). \tag{7.1}$$

Suppose there exists a constant $\delta > 0$ such that for $r = 2 + \delta$, the following conditions hold:

1. $\int \big|g(\boldsymbol{X}_1, \ldots, \boldsymbol{X}_m)\big|^r dF(\boldsymbol{X}_1) \ldots dF(\boldsymbol{X}_m) \le M_0 < \infty$ for some constant $M_0$;
2. $\mathbb{E}\big|g(\boldsymbol{X}_1, \ldots, \boldsymbol{X}_m)\big|^r \le M_1$ for some constant $M_1$;
3. $\{\boldsymbol{X}_t\}_{t\in\mathbb{Z}}$ is $\beta$-mixing with $\beta(n) = O\{n^{-(2+\delta')/\delta'}\}$ for some $0 < \delta' < \delta$.

Assuming that the above conditions hold, we then have

$$\frac{\sqrt{T}\{U_T(g) - \theta\}}{\sigma} \xrightarrow{\mathsf{d}} Z, \quad \text{as } T \to \infty,$$

where $Z \sim N(0, 1)$ is a standard Gaussian random variable.

**Lemma 7.3** (Yoshihara (1976))**.** Let $\{\boldsymbol{X}_t\}_{t\in\mathbb{Z}}$ be a $d$-dimensional stationary process with the marginal distribution function $F$, and $\boldsymbol{X}_1, \ldots, \boldsymbol{X}_T$ be a sequence of observations. Suppose $h(\cdot) : \mathbb{R}^d \times \mathbb{R}^d \to \mathbb{R}$ is a kernel function such that for some constants $\zeta > 0$ and $H > 0$, we have

$$\int\int |h(\boldsymbol{X}_1, \boldsymbol{X}_2)|^{2+\zeta} dF(\boldsymbol{X}_1) dF(\boldsymbol{X}_2) \le H, \tag{7.2}$$

$$\int |h(\boldsymbol{X}_1, \boldsymbol{X}_{1+k})|^{2+\zeta} d\mathbb{P}(\boldsymbol{X}_1, \boldsymbol{X}_{1+k}) \le H, \quad \text{for all} \;\; k \ge 0, \;\; k \in \mathbb{Z}, \tag{7.3}$$

where $\mathbb{P}(\boldsymbol{X}_{t_1}, \boldsymbol{X}_{t_2})$ is the joint distribution function of $(\boldsymbol{X}_{t_1}, \boldsymbol{X}_{t_2})$. For arbitrary random vectors $\{\boldsymbol{X}, \boldsymbol{Y}\}$, we define

$$h_1(\boldsymbol{X}) = \int h(\boldsymbol{X}, \boldsymbol{Y}) dF(\boldsymbol{Y}) - \int\int h(\boldsymbol{X}, \boldsymbol{Y}) dF(\boldsymbol{X}) dF(\boldsymbol{Y}),$$

$$h_2(\boldsymbol{X}, \boldsymbol{Y}) = h(\boldsymbol{X}, \boldsymbol{y}) - h_1(\boldsymbol{X}) - h_1(\boldsymbol{Y}) - \int\int h(\boldsymbol{X}, \boldsymbol{Y}) dF(\boldsymbol{X}) dF(\boldsymbol{Y}).$$

If the process $\{\boldsymbol{X}_t\}_{t\in\mathbb{Z}}$ is $\beta$-mixing with mixing coefficient $\beta(n) = O\{n^{-(2+\zeta')/\zeta'}\}$ for a constant $\zeta' \in (0, \zeta)$, then, for the $U$-statistic

$$U_T(h_2) := \frac{2}{T(T-1)} \sum_{t_1<t_2} h_2(\boldsymbol{X}_{t_1}, \boldsymbol{X}_{t_2}),$$

we have

$$\begin{aligned}\mathbb{E}\{TU_T(h_2)^2\} \leq & \frac{4}{T(T-1)^2} \sum_{1\leq t_1<t_2\leq T} \sum_{1\leq t_3<t_4\leq T} \Big|\mathbb{E}\{h_2(\boldsymbol{X}_{t_1}, \boldsymbol{X}_{t_2})h_2(\boldsymbol{X}_{t_3}, \boldsymbol{X}_{t_4})\}\Big| \\ \leq & \frac{4}{n^3} \sum_{t_1,t_2,t_3,t_4=1}^{T} \Big|\mathbb{E}\{h_2(\boldsymbol{X}_{t_1}, \boldsymbol{X}_{t_2})h_2(\boldsymbol{X}_{t_3}, \boldsymbol{X}_{t_4})\}\Big| = O(T^{-\lambda}),\end{aligned}$$

where $\lambda := \min\big(2(\zeta - \zeta')/\{\zeta'(2+\zeta)\}, 1\big)$.

**Lemma 7.4.** Let $\{\boldsymbol{X}_t\}_{t\in\mathbb{Z}}$ be a $d$-dimensional stationary process with the marginal distribution function $F$, $\boldsymbol{X}_1, \ldots, \boldsymbol{X}_T$ be a sequence of observations, and $\boldsymbol{X}_1^*, \ldots, \boldsymbol{X}_T^*$ be a block bootstrapped sample with block length $l \asymp T^{1-\epsilon_0}$ defined in Section 2.1. For a kernel function $h : \mathbb{R}^d \times \mathbb{R}^d \to \mathbb{R}$, define

$$U_T(h) = \frac{2}{T(T-1)} \sum_{t_1<t_2} h(\boldsymbol{X}_{t_1}, \boldsymbol{X}_{t_2}) \quad \text{and} \quad U_T^*(h) = \frac{2}{T(T-1)} \sum_{t_1<t_2} h(\boldsymbol{X}_{t_1}^*, \boldsymbol{X}_{t_2}^*)$$

to be the $U$-statistics based on the observed sample and bootstrap sample, respectively. Now supposing that $h$ satisfies (7.2) and (7.3), and the process $\{\boldsymbol{X}_t\}_{t\in\mathbb{Z}}$ is $\beta$-mixing with mixing coefficient $\beta(n) = O\{n^{-(2+\zeta')/\zeta'}\}$ for a constant $\zeta' \in (0, \zeta)$, we have

$$\Big|\mathrm{Var}^*\left\{\sqrt{T}U_T^*(h)\right\} - \mathrm{Var}\left\{\sqrt{T}U_T(h)\right\}\Big| = o_P(1),$$

where $\mathrm{Var}^*$ is the variance operator of the resampling distribution $\mathbb{P}^*$ conditional on $\boldsymbol{X}_1, \ldots, \boldsymbol{X}_T$.

*Proof.* We define $\omega := \displaystyle\int\int h(\boldsymbol{X}, \boldsymbol{Y})dF(\boldsymbol{X})dF(\boldsymbol{Y})$. Using Hoeffding's decomposition, we have

$$U_T^*(h) = \omega + \frac{2}{T}\sum_{t=1}^{T} h_1(\boldsymbol{X}_t^*) + U_T^*(h_2).$$

The fact that for two random variables $X$ and $Y$, we have $\mathrm{Var}(X+Y) = \mathrm{Var}(X) + \mathrm{Var}(Y) + 2\,\mathrm{Cov}(X, Y) \leq \mathrm{Var}(X) + \mathrm{Var}(Y) + 2\sqrt{\mathrm{Var}(X)}\sqrt{\mathrm{Var}(Y)}$, yields

$$\begin{aligned}\mathrm{Var}^*\{\sqrt{T}U_T^*(h)\} \leq & \mathrm{Var}^*\Big\{\frac{2}{\sqrt{T}}\sum_{t=1}^{T} h_1(\boldsymbol{X}_t^*)\Big\} + \mathrm{Var}^*\Big\{\sqrt{T}U_T^*(h_2)\Big\} \\ & + 2\sqrt{\mathrm{Var}^*\Big\{\frac{2}{\sqrt{T}}\sum_{t=1}^{T} h_1(\boldsymbol{X}_t^*)\Big\}}\sqrt{\mathrm{Var}^*\Big\{\sqrt{T}U_T^*(h_2)\Big\}}. \end{aligned} \tag{7.4}$$

Similarly, using the fact that $\text{Var}(X+Y) = \text{Var}(X) + \text{Var}(Y) + 2\,\text{Cov}(X,Y) \geq \text{Var}(X) + \text{Var}(Y) - 2\sqrt{\text{Var}(X)}\sqrt{\text{Var}(Y)}$, we have

$$\begin{aligned}\text{Var}^*\{\sqrt{T}U_T^*(h)\} \geq \text{Var}^*\Big\{\frac{2}{\sqrt{T}}\sum_{t=1}^T h_1(\boldsymbol{X}_t^*)\Big\} + \text{Var}^*\Big\{\sqrt{T}U_T^*(h_2)\Big\} \\ - 2\sqrt{\text{Var}^*\Big\{\frac{2}{\sqrt{T}}\sum_{t=1}^T h_1(\boldsymbol{X}_t^*)\Big\}}\sqrt{\text{Var}^*\Big\{\sqrt{T}U_T^*(h_2)\Big\}}.\end{aligned} \tag{7.5}$$

By Theorem 2.3 of Shao and Yu (1993), regarding $h_1$, we have

$$\Big|\text{Var}^*\Big\{\frac{2}{\sqrt{T}}\sum_{t=1}^T h_1(\boldsymbol{X}_t^*)\Big\} - \text{Var}\Big\{\frac{2}{\sqrt{T}}\sum_{t=1}^T h_1(\boldsymbol{X}_t)\Big\}\Big| \xrightarrow{\text{a.s.}} 0 \tag{7.6}$$

On the other hand, by Lemma 7.3, we have $\text{Var}\{\sqrt{T}U_T(h_2)\} = o(1)$ and $\text{Var}^*\{\sqrt{T}U_T^*(h_2)\} = o_P(1)$. Combining them with (7.4) and (7.5), we have

$$\text{Var}^*\{\sqrt{T}U_T^*(h)\} = \text{Var}^*\Big\{\frac{2}{\sqrt{T}}\sum_{t=1}^T h_1(\boldsymbol{X}_t^*)\Big\} + o_P(1). \tag{7.7}$$

Similar arguments yield that

$$\text{Var}\{\sqrt{T}U_T(h)\} = \text{Var}\Big\{\frac{2}{\sqrt{T}}\sum_{t=1}^T h_1(\boldsymbol{X}_t)\Big\} + o(1). \tag{7.8}$$

Combining (7.7) and (7.8), we obtain

$$\text{Var}^*\{\sqrt{T}U_T^*(h)\} - \text{Var}\{\sqrt{T}U_T(h)\} = \text{Var}^*\Big\{\frac{2}{\sqrt{T}}\sum_{t=1}^T h_1(\boldsymbol{X}_t^*)\Big\} - \text{Var}\Big\{\frac{2}{\sqrt{T}}\sum_{t=1}^T h_1(\boldsymbol{X}_t)\Big\} + o_P(1).$$

Combining the above equation with (7.6) completes the proof. □

**Lemma 7.5.** Let $\{\boldsymbol{X}_t\}_{t\in\mathbb{Z}}$ be a stationary sequence of $\phi$-mixing random vectors. Suppose the $\phi$-mixing coefficient satisfies Assumption **(A3)**. Then we have

$$\|\mathbb{E}\widehat{\mathbf{T}} - \mathbf{T}\|_{\max} = O\Big(1/T\Big),$$

where $\widehat{\mathbf{T}}$ and $\mathbf{T}$ are sample and population Kendall's tau matrix defined in (2.1).

*Proof.* For any two constant $1 \leq s < t \leq T$, we have

$$\mathbb{P}(X_{tj} - X_{sj} > 0, X_{tk} - X_{tk} > 0) = \mathbb{P}(X_{tj} > X_{sj}, X_{tk} > X_{sk}).$$

Let

$$-\infty = a_0 < -M < a_1 < \ldots < a_{h-1} < M < a_h = \infty$$

and

$$-\infty = b_0 < -M < b_1 < \ldots < b_{h-1} < M < b_h = \infty$$

be two pre-determined real sequences. Note that for $i_0 = 1, \ldots, h$, given $\{X_{sj} \in [a_{i_0-1}, a_{i_0}]\}$, the event $\{X_{tj} > X_{sj}\}$ implies the event $\{X_{tj} > a_{i_0-1}\}$. This yields

$$\mathbb{P}(X_{tj} > X_{sj}, X_{tk} > X_{sk}) \leq \sum_{i_0,j_0} \mathbb{P}(X_{tj} > a_{i_0-1}, X_{tk} > b_{j_0-1} \mid X_{sj} \in [a_{i_0-1}, a_{i_0}], X_{sk} \in [b_{j_0-1}, b_{j_0}]) \cdot$$
$$\mathbb{P}(X_{sj} \in [a_{i_0-1}, a_{i_0}], X_{sk} \in [b_{j_0-1}, b_{j_0}]).$$

On the other hand, given $\{X_{sj} \in [a_{i_0-1}, a_{i_0}]\}$, the event $\{X_{tj} > a_{i_0}\}$ implies the event $\{X_{tj} > X_{sj}\}$. Thus, we have

$$\mathbb{P}(X_{tj} > X_{sj}, X_{tk} > X_{sk}) \geq \sum_{i_0,j_0} \mathbb{P}(X_{tj} > a_{i_0}, X_{tk} > b_{j_0} \mid X_{sj} \in [a_{i_0-1}, a_{i_0}], X_{sk} \in [b_{j_0-1}, b_{j_0}]) \cdot$$
$$\mathbb{P}(X_{sj} \in [a_{i_0-1}, a_{i_0}], X_{sk} \in [b_{j_0-1}, b_{j_0}]).$$

Now, we define $\psi_h^U$ to be

$$\psi_h^U := \sum_{i_0,j_0} \mathbb{P}(X_{tj} > a_{i_0-1}, X_{tk} > b_{i_0-1}) \mathbb{P}(X_{sj} \in [a_{i_0-1}, a_{i_0}], X_{sk} \in [b_{j_0-1}, b_{j_0}]),$$

and similarly define $\psi_h^L$ to be

$$\psi_h^L = \sum_{i_0,j_0} \mathbb{P}(X_{tj} > a_{i_0}, X_{tk} > b_{i_0}) \mathbb{P}(X_{sj} \in [a_{i_0-1}, a_{i_0}], X_{sk} \in [b_{j_0-1}, b_{j_0}]).$$

Let $\psi_h$ be either $\psi_h^U$ or $\psi_h^L$ with regard to the sign of $\mathbb{P}(X_{tj} > X_{sj}, X_{tk} > X_{sk}) - \psi_h^L$:

$$\psi_h = \begin{cases} \psi_h^L, & \text{if } \mathbb{P}(X_{tj} > X_{sj}, X_{tk} > X_{sk}) > \psi_h^L; \\ \psi_h^U, & \text{Otherwise.} \end{cases}.$$

Without loss of generality, supposing that we have $\mathbb{P}(X_{tj} > X_{sj}, X_{tk} > X_{sk}) > \psi_h^L$, it follows that

$$\Big| \mathbb{P}(X_{tj} > X_{sj}, X_{tk} > X_{sk}) - \psi_h \Big| = \mathbb{P}(X_{tj} > X_{sj}, X_{tk} > X_{sk}) - \psi_h^L$$
$$\leq \sum_{i_0,j_0} \Big| \mathbb{P}(X_{tj} > a_{j_0-1}, X_{tk} > b_{j_0-1} \mid X_{sj} \in [a_{i_0-1}, a_{i_0}], X_{sk} \in [b_{j_0-1}, b_{j_0}]) - \mathbb{P}(X_{tj} > a_{j_0}, X_{tk} > b_{j_0}) \Big| \cdot$$
$$\mathbb{P}(X_{sj} \in [a_{i_0-1}, a_{i_0}], X_{sk} \in [b_{j_0-1}, b_{j_0}])$$
$$\leq \phi(t-s) + \max_{i_0,j_0} |\mathbb{P}(X_{tj} > a_{j_0-1}, X_{tk} > b_{j_0-1}) - \mathbb{P}(X_{tj} > a_{j_0}, X_{tk} > b_{j_0})|.$$

Now let $h \to \infty$, $\max_{i=2}^{h-1} |a_i - a_{i-1}| \to 0$, $\max_{i=2}^{h-1} |b_i - b_{i-1}| \to 0$, and $M \to \infty$. By the definition of $\phi$-mixing coefficient, we have

$$\Big|\mathbb{P}(X_{tj} > X_{sj}, X_{tk} > X_{sk}) - \int \mathbb{P}(X_{tj} > a, X_{tk} > b) d\mathbb{P}(X_{sj} = a, X_{sk} = b)\Big| \leq \phi(s - t). \quad (7.9)$$

Moreover, letting $\boldsymbol{X}' = (X'_1, \ldots, X'_d)^T$ have the same distribution as $\boldsymbol{X}_1$ and independent of $(\boldsymbol{X}_s, \boldsymbol{X}_t)$, we have

$$d\mathbb{P}(X_{sj} = a, X_{sk} = b) = d\mathbb{P}(X'_j = a, X'_k = b).$$

This yields

$$\int \mathbb{P}(X_{tj} > a, X_{tk} > b) d\mathbb{P}(X_{sj} = a, X_{sk} = b) = \int \mathbb{P}(X_{tj} > a, X_{tk} > b) d\mathbb{P}(X'_j = a, X'_k = b).$$

Plugging the above equation into (7.9), we obtain

$$\Big|\mathbb{P}(X_{tj} > X_{sj}, X_{tk} > X_{sk}) - \int \mathbb{P}(X_{tj} > a, X_{tk} > b) d\mathbb{P}(X'_j = a, X'_k = b)\Big| \leq \phi(t - s).$$

Note that by the definition of conditional probability, we have

$$\int \mathbb{P}(X_{tj} > a, X_{tk} > b) d\mathbb{P}(X'_j = a, X'_k = b) = \mathbb{P}(X_{tj} - X'_j > 0, X_{tk} - X'_k > 0).$$

Thus, combining the above two equations, we have

$$\Big|\mathbb{P}(X_{tj} - X_{sj} > 0, X_{tk} - X_{sk} > 0) - \mathbb{P}(X_{tj} - X'_j > 0, X_{tk} - X'_k > 0)\Big| \leq \phi(t - s). \quad (7.10)$$

Using similar arguments, we can prove

$$\Big|\mathbb{P}(X_{tj} - X_{sj} < 0, X_{tk} - X_{sk} < 0) - \mathbb{P}(X_{tj} - X'_j < 0, X_{tk} - X'_k < 0)\Big| \leq \phi(t - s), \quad (7.11)$$
$$\Big|\mathbb{P}(X_{tj} - X_{sj} < 0, X_{tk} - X_{sk} > 0) - \mathbb{P}(X_{tj} - X'_j < 0, X_{tk} - X'_k > 0)\Big| \leq \phi(t - s), \quad (7.12)$$
$$\Big|\mathbb{P}(X_{tj} - X_{sj} > 0, X_{tk} - X_{sk} < 0) - \mathbb{P}(X_{tj} - X'_j > 0, X_{tk} - X'_k < 0)\Big| \leq \phi(t - s). \quad (7.13)$$

By definition, we have $\tau_{jk} = \mathbb{E}\{\text{sign}(X_{tj} - X'_j)(X_{tk} - X'_k)\}$. Applying the definition of expectation, we have

$$\begin{aligned}\tau_{jk} =& \mathbb{P}(X_{tj} - X'_j > 0, X_{tk} - X'_k > 0) + \mathbb{P}(X_{tj} - X'_j < 0, X_{tk} - X'_k < 0) - \\ & \mathbb{P}(X_{tj} - X'_j > 0, X_{tk} - X'_k < 0) - \mathbb{P}(X_{tj} - X'_j < 0, X_{tk} - X'_k > 0). \quad (7.14)\end{aligned}$$

By the same reason, we have

$$
\begin{aligned}
&\mathbb{E}\{\text{sign}(X_{tj}-X_{sj})(X_{tk}-X_{sk})\}\\
=&\mathbb{P}(X_{tj}-X_{sj}>0, X_{tk}-X'_{sk}>0)+\mathbb{P}(X_{tj}-X_{sj}<0, X_{tk}-X_{sk}<0)-\\
&\mathbb{P}(X_{tj}-X_{sj}>0, X_{tk}-X_{sk}<0)-\mathbb{P}(X_{tj}-X_{sj}<0, X_{tk}-X_{sk}>0).
\end{aligned} \tag{7.15}
$$

Now, by the definition of $\widehat{\tau}_{jk}$, we have

$$
\begin{aligned}
\Big|\mathbb{E}\widehat{\tau}_{jk}-\tau_{jk}\Big| =&\Big|\mathbb{E}\Big\{\frac{2}{T(T-1)}\sum_{s<t}\text{sign}(X_{tj}-X_{sj})(X_{tk}-X_{sk})\Big\}-\tau_{jk}\Big|\\
\leq&\frac{2}{T(T-1)}\sum_{s<t}\Big|\mathbb{E}\text{sign}(X_{tj}-X_{sj})(X_{tk}-X_{sk})-\tau_{jk}\Big|.
\end{aligned}
$$

Plugging (7.14) and (7.15) into the above equation, and applying (7.10) - (7.13), we obtain

$$
\Big|\mathbb{E}\widehat{\tau}_{jk}-\tau_{jk}\Big|\leq\frac{2}{T(T-1)}\sum_{s<t}\{4\phi(t-s)\}=\frac{8\sum_{t=1}^{T}(T-t)\phi(t)}{T(T-1)}=O\Big(\frac{1}{T}\Big). \tag{7.16}
$$

The last inequality is because by Assumption **(A3)**, we have

$$
\sum_{t=1}^{T}(T-t)\phi(t)\leq\sum_{t=1}^{T}\frac{T-t}{t^{1+\epsilon}}\leq T\sum_{t=1}^{\infty}\frac{1}{t^{1+\epsilon}}=O(T).
$$

This completes the proof. □

**Lemma 7.6.** Let $\{\boldsymbol{X}_t\}_{t\in\mathbb{Z}}$ be a stationary sequence of $\phi$-mixing random vectors. Suppose the $\phi$-mixing coefficient satisfies Assumption **(A3)**. Then we have

$$
\|\widehat{\mathbf{T}}-\mathbf{T}\|_{\max}=O_P\Big(\sqrt{\frac{\log d}{T}}\Big),
$$

where $\widehat{\mathbf{T}}$ and $\mathbf{T}$ are the sample and population Kendall's tau matrix based on $\{\boldsymbol{X}_t\}_{t=1}^{T}$.

*Proof.* Consider the following function

$$
f_{jk}(\boldsymbol{X}_1,\ldots,\boldsymbol{X}_T):=\frac{2}{T-1}\sum_{t<t'}\text{sign}(X_{tj}-X_{t'j})\text{sign}(X_{tk}-X_{t'k})=T\cdot\widehat{\tau}_{jk}.
$$

We have

$$
\begin{aligned}
&\Big|f_{jk}(\boldsymbol{X}_1,\ldots,\boldsymbol{X}_i,\ldots,\boldsymbol{X}_T)-f_{jk}(\boldsymbol{X}_1,\ldots,\boldsymbol{X}'_i,\ldots,\boldsymbol{X}_T)\Big|\\
=&\frac{2}{T-1}\Big|\sum_{t\neq i}\text{sign}(X_{ij}-X_{tj})\text{sign}(X_{ik}-X_{tk})-\sum_{t\neq i'}\text{sign}(X_{i'j}-X_{tj})\text{sign}(X_{i'k}-X_{tk})\Big|\\
\leq&\frac{2}{T-1}\{2(T-1)\}=4.
\end{aligned}
$$

Thus, $f_{jk}$ is $c$-Lipschitz with respect to the Hamming metric. By Lemma 7.1, we have

$$\mathbb{P}\Big(T|\widehat{\tau}_{jk} - \mathbb{E}\widehat{\tau}_{jk}| \geq \epsilon\Big) \leq 2\exp\Big[-\frac{\epsilon^2}{8T\{1+2\sum_{l=1}^{\infty}\phi(l)\}}\Big],$$

for any $\epsilon > 0$. Here $\sum_{l=1}^{\infty}\phi(l) < \infty$ is guaranteed by Assumption **(A3)**. Thus, we have

$$\mathbb{P}\Big(\|\widehat{\mathbf{T}} - \mathbb{E}\widehat{\mathbf{T}}\|_{\max} \geq \epsilon\Big) \leq \sum_{j,k=1}^{d}\mathbb{P}\Big(|\widehat{\tau}_{jk} - \mathbb{E}\widehat{\tau}_{jk}| \geq \epsilon\Big) \leq 2\exp\Big[2\log d - \frac{T\epsilon^2}{8\{1+2\sum_{l=1}^{\infty}\phi(l)\}}\Big].$$

Setting $\epsilon = \sqrt{[24\{1+2\sum_{l=1}^{\infty}\phi(l)\}\log d]/T}$, we have

$$\|\widehat{\mathbf{T}} - \mathbb{E}\widehat{\mathbf{T}}\|_{\max} = O_P\Big(\sqrt{\frac{\log d}{T}}\Big).$$

Combining the above equation with Lemma 7.5 completes the proof. □

**Lemma 7.7.** [Theorem 1 in Doukhan and Neumann (2007)] Suppose that $X_1, \ldots, X_T$ are real-valued random variables with mean 0, defined on a common probability space $(\Omega, \mathcal{A}, \mathbb{P})$. Let $\Psi : \mathbb{N}^2 \to \mathbb{N}$ be one of the following functions:

$(a)$. $\Psi(u, v) = 2v$,

$(b)$. $\Psi(u, v) = u + v$,

$(c)$. $\Psi(u, v) = uv$,

$(d)$. $\Psi(u, v) = \alpha(u + v) + (1-\alpha)uv$, for some $\alpha \in (0, 1)$.

We assume that there exist constants $K, M, L_1, L_2 > 0$, $a, b \geq 0$, and a non-increasing sequence of real coefficients $\{\rho(n)\}_{n\geq 0}$ such that for any $u$-tuple $(s_1, \ldots, s_u)$ and $v$-tuple $(t_1, \ldots, t_v)$ with $1 \leq s_1 \leq \cdots \leq s_u \leq t_1 \leq \cdots \leq t_v \leq T$, the following inequalities hold:

$$\left|\mathrm{Cov}\left(\prod_{i=1}^{u} X_{s_i}, \prod_{j=1}^{v} X_{t_j}\right)\right| \leq K^2 M^{u+v}\{(u+v)!\}^b \Psi(u, v)\rho(t_1 - s_u), \tag{7.17}$$

where for the sequence $\{\rho(n)\}_{n\geq 0}$, we require that

$$\sum_{s=0}^{\infty}(s+1)^k \rho(s) \leq L_1 L_2^k (k!)^a, \quad \forall k \geq 0. \tag{7.18}$$

We also assume that the following moment condition holds:

$$\mathbb{E}|X_t|^k \leq (k!)^b M^k, \text{ for all } t = 1, \ldots, T. \tag{7.19}$$

Let $S_T = \sum_{t=1}^T X_t$. Then, for all $x > 0$, we have

$$\mathbb{P}(S_T \geq x) \leq \exp\left\{-\frac{x^2}{C_1 T + C_2 x^{(2a+2b+3)/(a+b+2)}}\right\}, \tag{7.20}$$

where $C_1$ and $C_2$ are constants depending on $K, M, L_1, L_2, a$, and $b$:

$$C_1 = 2^{a+b+3} K^2 M^2 L_1 (K^2 \vee 2), \quad C_2 = 2\{M L_2 (K^2 \vee 2)\}^{1/(a+b+2)}. \tag{7.21}$$

**Lemma 7.8.** Let $\{\boldsymbol{X}_t\}_{t\in\mathbb{Z}}$ be a $d$ dimensional stationary $\phi$-mixing process satisfying Assumptions **(A6)**, **(A7)**, and **(A9)**. Let $\widehat{\mathbf{R}} = \mathrm{diag}(\widehat{\sigma}_{\mathrm{M},1}, \ldots, \widehat{\sigma}_{\mathrm{M},d})$ be a diagonal matrix of sample median absolute deviations based on $\{\boldsymbol{X}_t\}_{t=1}^T$, and $\mathbf{R} = \mathrm{diag}\{\sigma_{\mathrm{M}}(X_{11}), \ldots, \sigma_{\mathrm{M}}(X_{1d})\}$ be its population counterpart. Then we have

$$\|\widehat{\mathbf{R}} - \mathbf{R}\|_{\max} = O_P\Big(\sqrt{\frac{\log d}{T}}\Big).$$

*Proof.* We first focus on a marginal process $\{X_{tj}\}_{t=1}^T$. For notational brevity, we suppress the index $j$ and denote the process as $\{X_t\}_{t=1}^T$. Define $X = X_1$. Let $F$ be the distribution function of $X$ and $F_T$ be the empirical distribution of $\{X_t\}_{t=1}^T$ and $F_T^{-1}(q) := \widehat{Q}(\{X_t\}; q)$ for any $q \in [0,1]$. By the definition of $\widehat{Q}(\cdot)$ in (2.3), we have, for any $\epsilon \in [0,1]$,

$$\epsilon \leq F_T\{F_T^{-1}(\epsilon)\} \leq \epsilon + \frac{1}{T}.$$

This implies that

$$\mathbb{P}\Big\{\widehat{Q}(\{X_t\}; q) - Q(X; q) \geq u\Big\} = \mathbb{P}\Big\{F_T^{-1}(q) - F^{-1}(q) \geq u\Big\} \leq \mathbb{P}\Big[q + \frac{1}{T} \geq F_T\{u + F^{-1}(q)\}\Big].$$

By the definition of $F_T$, we further have

$$\begin{aligned}
&\mathbb{P}\Big\{\widehat{Q}(\{X_t\}; q) - Q(X; q) \geq u\Big\} \leq \mathbb{P}\Big[\sum_{t=1}^T I\{X_t \leq F^{-1}(q) + u\} \leq Tq + 1\Big] \\
=&\mathbb{P}\Big(\sum_{t=1}^T \Big[-I\{X_t \leq F^{-1}(q) + u\} + F\{F^{-1}(q) + u\}\Big] \geq T\Big[F\{F^{-1}(q) + u\} - q - \frac{1}{T}\Big]\Big).
\end{aligned}$$

Since $\{\boldsymbol{X}_t\}_{t\in\mathbb{Z}}$ is $\phi$-mixing, the process $\{-I\{X_t \leq F^{-1}(q) + u\} + F\{F^{-1}(q) + u\}\}_{t\in\mathbb{Z}}$ is also $\phi$-mixing. By Lemma 6 in Doukhan and Louhichi (1999), $\{-I\{X_t \leq F^{-1}(q) + u\} + F\{F^{-1}(q) + u\}\}_{t\in\mathbb{Z}}$ satisfies (7.17) with $K = 2$, $M = 1$, $b = 0$, any of the four $\Psi$ functions, and

$$\rho(n) = \phi(n) \leq C_1 \exp(-C_2 n^r).$$

By Proposition 8 in Doukhan and Neumann (2007), (7.18) is satisfied with $a = \max(1, 1/r)$ and some constants $L_1$ and $L_2$. Since $-I\{X_t \le F^{-1}(q) + u\} + F\{F^{-1}(q) + u\}$ is bounded, (7.19) is also satisfied with $b = 0$. Thus, applying Lemma 7.7, we have

$$\mathbb{P}\Big\{\widehat{Q}(\{X_t\}; q) - Q(X; q) \ge u\Big\} \le \exp\Big(-\psi\Big(F\{F^{-1}(q) + u\} - q - \frac{1}{T}\Big)\Big), \tag{7.22}$$

for $F\{F^{-1}(q) + u\} - q - 1/T > 0$, where

$$\psi(x) := \frac{Tx^2}{C_1 + C_2 T^{(a+1)/(a+2)} x^{(2a+3)/(a+2)}},$$

for $x > 0$, $a = \max(1, 1/r)$, and some absolute constants $C_1$ and $C_2$. On the other hand, we have

$$\mathbb{P}\Big\{\widehat{Q}(\{X_t\}; q) - Q(X; q) \le -u\Big\} = \mathbb{P}\Big\{F_T^{-1}(q) - F^{-1}(q) \le -u\Big\} \le \mathbb{P}\Big[q \le F_T\{F^{-1}(q) - u\}\Big]$$
$$= \mathbb{P}\Big(\sum_{t=1}^{T} \Big[I\{X_t \le F^{-1}(q) - u\} - F\{F^{-1}(q) - u\}\Big] \ge T\Big[q - F\{F^{-1}(q) - u\}\Big]\Big).$$

By similar arguments, we have

$$\mathbb{P}\Big\{\widehat{Q}(\{X_t\}; q) - Q(X; q) \le -u\Big\} \le \exp\Big(-\psi\Big(q - F\{F^{-1}(q) - u\}\Big)\Big). \tag{7.23}$$

Combining (7.22) and (7.23), we have

$$\mathbb{P}\Big\{\Big|\widehat{Q}(\{X_t\}; q) - Q(X; q)\Big| \ge u\Big\}$$
$$\le \exp\Big(-\psi\Big(F\{F^{-1}(q) + u\} - q - \frac{1}{T}\Big)\Big) + \exp\Big(-\psi\Big(q - F\{F^{-1}(q) - u\}\Big)\Big), \tag{7.24}$$

for $F\{F^{-1}(q) + u\} - q - 1/T > 0$.

Next, we continue to derive exponential tail probabilities for $\widehat{\sigma}_{\mathrm{M}}(\{X_t\}_{t=1}^T)$. We write $\widehat{m} := \widehat{Q}(\{X_t\}_{t=1}^T; 1/2)$ and $m := Q(X_1; 1/2)$ to be the sample and population medians. Let $F_1$ and $F_2$ be the distribution functions of $X$ and $|X - Q(X; 1/2)|$. By the definition of $\widehat{\sigma}_{\mathrm{M}}$, we have

$$\mathbb{P}\Big\{\widehat{\sigma}_{\mathrm{M}}\Big(\{X_t\}_{t=1}^T\Big) - \sigma_{\mathrm{M}}(X) > u\Big\} = \mathbb{P}\Big\{\widehat{Q}\Big(\Big\{|X_t - \widehat{m}|\Big\}_{t=1}^T; \frac{1}{2}\Big) - Q\Big(|X - m|; \frac{1}{2}\Big) > u\Big\}$$
$$\le \mathbb{P}\Big\{\widehat{Q}\Big(\Big\{|X_t - m|\Big\}_{t=1}^T; \frac{1}{2}\Big) + |\widehat{m} - m| - Q\Big(|X - m|; \frac{1}{2}\Big) > u\Big\}$$
$$\le \mathbb{P}\Big\{\widehat{Q}\Big(\Big\{|X_t - m|\Big\}_{t=1}^T; \frac{1}{2}\Big) - Q\Big(|X - m|; \frac{1}{2}\Big) > \frac{u}{2}\Big\} + \mathbb{P}\Big(|\widehat{m} - m| > \frac{u}{2}\Big). \tag{7.25}$$

On the other hand, using the same technique, we have

$$\begin{aligned}
&\mathbb{P}\Big\{\widehat{\sigma}_{\mathrm{M}}\Big(\{X_t\}_{t=1}^T\Big)-\sigma_{\mathrm{M}}(X)<-u\Big\}=\mathbb{P}\Big\{\widehat{Q}\Big(\Big\{|X_t-\widehat{m}|\Big\}_{t=1}^T;\frac{1}{2}\Big)-Q\Big(|X-m|;\frac{1}{2}\Big)<-u\Big\}\\
\le&\mathbb{P}\Big\{\widehat{Q}\Big(\Big\{|X_t-m|\Big\}_{t=1}^T;\frac{1}{2}\Big)-|\widehat{m}-m|-Q\Big(|X-m|;\frac{1}{2}\Big)<-u\Big\}\\
\le&\mathbb{P}\Big\{\widehat{Q}\Big(\Big\{|X_t-m|\Big\}_{t=1}^T;\frac{1}{2}\Big)-Q\Big(|X-m|;\frac{1}{2}\Big)<-\frac{u}{2}\Big\}+\mathbb{P}\Big(|\widehat{m}-m|>\frac{u}{2}\Big). \qquad (7.26)
\end{aligned}$$

Combining (7.25) and (7.26), we have

$$\begin{aligned}
&\mathbb{P}\Big\{|\widehat{\sigma}_{\mathrm{M}}\Big(\{X_t\}_{t=1}^T\Big)-\sigma_{\mathrm{M}}(X)|>u\Big\}\\
\le&\mathbb{P}\Big\{\Big|\widehat{Q}\Big(\Big\{|X_t-m|\Big\}_{t=1}^T;\frac{1}{2}\Big)-Q\Big(|X-m|;\frac{1}{2}\Big)\Big|>\frac{u}{2}\Big\}+2\mathbb{P}\Big(|\widehat{m}-m|>\frac{u}{2}\Big). \qquad (7.27)
\end{aligned}$$

Now applying Inequality (7.24), we have

$$\begin{aligned}
&\mathbb{P}\Big\{|\widehat{Q}\Big(\Big\{|X_t-m|\Big\}_{t=1}^T;\frac{1}{2}\Big)-Q\Big(|X-m|;\frac{1}{2}\Big)|>\frac{u}{2}\Big\}\\
\le&\exp\Big(-\psi\Big[F_2\Big\{F_2^{-1}\Big(\frac{1}{2}\Big)+\frac{u}{2}\Big\}-\frac{1}{2}-\frac{1}{T}\Big]\Big)+\exp\Big(-\psi\Big[\frac{1}{2}-F_2\Big\{F_2^{-1}\Big(\frac{1}{2}\Big)-\frac{u}{2}\Big\}\Big]\Big)\\
\le&\exp\Big\{-\psi\Big(\frac{\eta u}{2}-\frac{1}{T}\Big)\Big\}+\exp\Big\{-\psi\Big(\frac{\eta u}{2}\Big)\Big\}, \qquad (7.28)
\end{aligned}$$

whenever $F_2\{F_2^{-1}(1/2)+u/2\}-1/2>1/T$. Here the last inequality is due to Assumption **(A9)** and the fact that $\psi$ is non-decreasing. Similarly, we also have

$$\begin{aligned}
&\mathbb{P}\Big(|\widehat{m}-m|>\frac{u}{2}\Big)\\
\le&\exp\Big(-\psi\Big[F_1\Big\{F_1^{-1}\Big(\frac{1}{2}\Big)+\frac{u}{2}\Big\}\Big]-\frac{1}{2}-\frac{1}{T}\Big)+\exp\Big(-\psi\Big[\frac{1}{2}-F_1\Big\{F_1^{-1}\Big(\frac{1}{2}\Big)-\frac{u}{2}\Big\}\Big]\Big)\\
\le&\exp\Big\{-\psi\Big(\frac{\eta u}{2}-\frac{1}{T}\Big)\Big\}+\exp\Big\{-\psi\Big(\frac{\eta u}{2}\Big)\Big\}, \qquad (7.29)
\end{aligned}$$

whenever $F_1\{F_1^{-1}(1/2)+u/2\}-1/2>1/T$. Again the last inequality is due to Assumption **(A9)** and the fact that $f$ is nondecreasing. Here we recall that $F_1$ and $F_2$ are the distribution functions of $X$ and $|X-Q(X;1/2)|$. Combining Inequalities (7.27), (7.28), and (7.29), we have

$$\begin{aligned}
\mathbb{P}\Big\{\Big|\widehat{\sigma}_{\mathrm{M}}\Big(\{X_t\}_{t=1}^T\Big)-\sigma_{\mathrm{M}}(X)\Big|>u\Big\}\le&3\exp\Big\{-\psi\Big(\frac{\eta u}{2}-\frac{1}{T}\Big)\Big\}+3\exp\Big\{-\psi\Big(\frac{\eta u}{2}\Big)\Big\}\\
\le&6\exp\Big\{-\psi\Big(\frac{\eta u}{2}-\frac{1}{T}\Big)\Big\},
\end{aligned}$$

whenever we have $0 < u/2 < \kappa$ and $\eta u/2 > 1/T$. Now we switch the focus back to the entire matrix $\widehat{\mathbf{R}}$. By the sub-additivity of probability measures, we have

$$\begin{aligned}\mathbb{P}\Big(\|\widehat{\mathbf{R}} - \mathbf{R}\|_{\max} > u\Big) \leq & \sum_{j=1}^{d} \mathbb{P}\Big\{\Big|\widehat{\sigma}_{\mathrm{M}}\Big(\{X_{tj}\}_{t=1}^{T}\Big) - \sigma_{\mathrm{M}}(X_{1j})\Big| > u\Big\} \\ \leq & 6\exp\Big\{2\log d - \psi\Big(\frac{\eta u}{2} - \frac{1}{T}\Big)\Big\}. \end{aligned} \tag{7.30}$$

We recall that by the definition of the function $\psi(\cdot)$, we have

$$\psi\Big(\frac{\eta u}{2} - \frac{1}{T}\Big) = \frac{T\big(\frac{\eta u}{2} - \frac{1}{T}\big)^2}{C_1 + C_2 T^{(a+1)/(a+2)}\big(\frac{\eta u}{2} - \frac{1}{T}\big)^{(2a+3)/(a+2)}}.$$

To simplify the denominator on the right-hand side of the above equation, we require that

$$C_1 \geq C_2 T^{(a+1)/(a+2)}\Big(\frac{\eta u}{2} - \frac{1}{T}\Big)^{(2a+3)/(a+2)}. \tag{7.31}$$

Then we have $\psi(\eta u/2 - 1/T) \geq T/(2C_1)(\eta u/2 - 1/T)^2$. Plugging this into (7.30), we obtain

$$\mathbb{P}\Big(\|\widehat{\mathbf{R}} - \mathbf{R}\|_{\max} > u\Big) \leq 6\exp\Big\{2\log d - \frac{T}{2C_1}\Big(\frac{\eta u}{2} - \frac{1}{T}\Big)^2\Big\}. \tag{7.32}$$

Next we select a proper $u$ to derive the rate of convergence. To this end, we set

$$2\log d - \frac{T}{2C_1}\Big(\frac{\eta u}{2} - \frac{1}{T}\Big)^2 = -\log d.$$

This leads to

$$u = \frac{2}{\eta}\Big(\sqrt{\frac{6C_1 \log d}{T}} + \frac{1}{T}\Big). \tag{7.33}$$

Plugging the above equation into (7.31), we get

$$C_1 \geq 6^{(2a+3)/(2a+4)} C_2 \Big\{\frac{(\log d)^{2a+3}}{T}\Big\}^{1/(2a+4)}.$$

Thus, (7.31) holds as long as we have $\log d = o[T^{1/(2a+3)}]$. By Assumption **(A7)**, (7.31) holds. Plugging (7.33) into (7.32), we get

$$\mathbb{P}\Big\{\|\widehat{\mathbf{R}} - \mathbf{R}\|_{\max} > \frac{2}{\eta}\Big(\sqrt{\frac{6C_1 \log d}{T}} + \frac{1}{T}\Big)\Big\} \leq \frac{6}{d}.$$

Thus, as $T$ and $d$ both go to infinity, we have

$$\|\widehat{\mathbf{R}} - \mathbf{R}\|_{\max} = O_P\Big(\sqrt{\frac{\log d}{T}}\Big).$$

This completes the proof. □

**Lemma 7.9.** Let $\{\boldsymbol{X}_t\}_{t=1}^T$ be a $d$ dimensional stationary process satisfying Assumptions **(A6) - (A9)**. We then have

$$\|\widehat{\mathbf{D}} - \mathbf{D}\|_{\max} = O_P\Big(\sqrt{\frac{\log d}{T}}\Big),$$

where $\widehat{\mathbf{D}}$ is defined in Equation (2.13).

*Proof.* Define $\widehat{\mathbf{R}} = \mathrm{diag}(\widehat{\sigma}_{\mathrm{M},1}, \ldots, \widehat{\sigma}_{\mathrm{M},d})$, $\mathbf{R} = \mathrm{diag}\{\sigma_{\mathrm{M}}(X_{11}), \ldots, \sigma_{\mathrm{M}}(X_{1d})\}$, $\widehat{c}^{\mathrm{M}} = \widehat{\sigma}_1/\widehat{\sigma}_{\mathrm{M},1}$, and $c^{\mathrm{M}} = \sqrt{\mathbf{\Sigma}_{11}}/\sigma_{\mathrm{M}}(X_{11})$. We have

$$\begin{aligned}\|\widehat{\mathbf{D}} - \mathbf{D}\|_{\max} &= \|\widehat{c}^{\mathrm{M}}\widehat{\mathbf{R}} - c^{\mathrm{M}}\mathbf{R}\|_{\max} \le \|\widehat{c}^{\mathrm{M}}(\widehat{\mathbf{R}} - \mathbf{R})\|_{\max} + \|(\widehat{c}^{\mathrm{M}} - c^{\mathrm{M}})\mathbf{R}\|_{\max} \\ &\le |\widehat{c}^{\mathrm{M}}|\|\widehat{\mathbf{R}} - \mathbf{R}\|_{\max} + C|\widehat{c}^{\mathrm{M}} - c^{\mathrm{M}}|. \end{aligned} \tag{7.34}$$

By Lemma 7.8, we have

$$\|\widehat{\mathbf{R}} - \mathbf{R}\|_{\max} = O_P\Big(\sqrt{\frac{\log d}{T}}\Big). \tag{7.35}$$

Thus, specifically, we have

$$\widehat{\sigma}_{\mathrm{M},1} \xrightarrow{P} \sigma_{\mathrm{M}}(X_{11}). \tag{7.36}$$

We can rewrite $\widehat{\sigma}_1^2$ as

$$\widehat{\sigma}_1^2 = \frac{1}{T-1}\sum_{t=1}^T (X_{t1} - \bar{X}_{T1})^2 = \frac{2}{T(T-1)}\sum_{t<t'} h(X_{t1}, X_{t'1}),$$

where $\bar{X}_{T1} := \sum_{t=1}^T X_{tj}/T$, and $h(X_{t1}, X_{t'1}) = (X_{t1} - X_{t'1})^2/2$. Thus, $\widehat{\sigma}_1^2$ is a $U$-statistic with kernel function $h$. Using Lemma 7.2 with Assumptions **(A6)** and **(A8)**, we have $\sqrt{T}(\widehat{\sigma}_1^2 - \mathbf{\Sigma}_{11}) \xrightarrow{\mathsf{d}} Z_1$ where $Z_1$ is a Gaussian random variable with mean 0. Using the delta method, we have $\sqrt{T}(\widehat{\sigma}_1 - \sqrt{\mathbf{\Sigma}_{11}}) \xrightarrow{\mathsf{d}} Z_2$ for another mean 0 Gaussian random variable $Z_2$. Combining this with (7.36) and applying Slutsky's theorem, we have $\sqrt{T}(\widehat{c}^{\mathrm{M}} - c^{\mathrm{M}}) \xrightarrow{\mathsf{d}} Z_3$ for some Gaussian random variable $Z_3$. Thus, we have

$$|\widehat{c}^{\mathrm{M}} - c^{\mathrm{M}}| = O_P\big(1/\sqrt{T}\big). \tag{7.37}$$

Combining (7.34), (7.35), and (7.37), we have the desired result. □

## 7.2 Proof of Theorem 3.1

*Proof.* Denote $\boldsymbol{a} = \mathbf{D}\mathbf{w}$. Using Taylor expansion entry-wise on $\sin(\pi\widehat{\mathbf{T}}/2)$ at $\sin(\pi\mathbf{T}/2)$, we have

$$\mathbf{w}^{\mathsf{T}}(\widehat{\boldsymbol{\Sigma}} - \boldsymbol{\Sigma})\mathbf{w} = \boldsymbol{a}^{\mathsf{T}}\Big\{\sin(\frac{\pi}{2}\widehat{\mathbf{T}}) - \sin(\frac{\pi}{2}\mathbf{T})\Big\}\boldsymbol{a}$$
$$= \underbrace{\boldsymbol{a}^{\mathsf{T}}\Big\{\cos(\frac{\pi}{2}\mathbf{T}) \circ \frac{\pi}{2}(\widehat{\mathbf{T}} - \mathbf{T})\Big\}\boldsymbol{a}}_{A_1} + \underbrace{\boldsymbol{a}^{\mathsf{T}}\Big\{-\frac{1}{2}[\sin(\theta_{jk})] \circ (\frac{\pi}{2})^2(\widehat{\mathbf{T}} - \mathbf{T}) \circ (\widehat{\mathbf{T}} - \mathbf{T})\Big\}\boldsymbol{a}}_{A_2},$$

where for each $j, k \in \{1, \ldots, d\}$, $\theta_{jk}$ lies between $\tau_{jk}$ and $\widehat{\tau}_{jk}$. Using Lemma 7.6 and assumption **(A4)**, we have

$$A_2 \le \frac{\pi^2}{8}\|\boldsymbol{a}\|_1^2\|\widehat{\mathbf{T}} - \mathbf{T}\|_{\max}^2 = O_P\Big(\frac{\log d}{T}\Big) = o_P\Big(\frac{\sigma}{\sqrt{T}}\Big). \tag{7.38}$$

Here the first inequality is due to the fact that for any vectors $\mathbf{v}_1, \mathbf{v}_2 \in \mathbb{R}^d$ and matrix $\mathbf{M} \in \mathbb{R}^{d\times d}$,

$$|\mathbf{v}_1^{\mathsf{T}}\mathbf{M}\mathbf{v}_2| \le \|\mathbf{v}_1\|_1\|\mathbf{M}\mathbf{v}_2\|_\infty \le \|\mathbf{M}\|_{\max}\|\mathbf{v}_1\|_1\|\mathbf{v}_2\|_1. \tag{7.39}$$

Next, we focus on $A_1$. We can expand $A_1$ by

$$A_1 = \underbrace{\frac{2}{T(T-1)}\sum_{t<t'} g(\boldsymbol{R}_t, \boldsymbol{R}_{t'})}_{U_T} - \boldsymbol{a}^{\mathsf{T}}\Big\{\cos(\frac{\pi}{2}\mathbf{T}) \circ \frac{\pi}{2}\mathbf{T}\Big\}\boldsymbol{a}, \tag{7.40}$$

where $g(\cdot)$ is defined in Equation (3.1). Note that $U_T$ is a $U$-statistic of order 2 and the kernel function $g(\cdot)$ satisfying

$$\Big|g(\boldsymbol{R}_t, \boldsymbol{R}_{t'})\Big| \le \frac{\pi}{2}\max_{jk}\Big|\text{sign}(R_{tj} - R_{t'j})\text{sign}(R_{tk} - R_{t'k})\Big|\|\boldsymbol{a}\|_1^2 \le \frac{\pi}{2}\|\mathbf{D}\|_{\max}\|\mathbf{w}\|_1^2 \le \frac{\pi}{2}.$$

Thus $g(\cdot)$ is a bounded kernel function. Assumption **(A3)** guarantees that $\{\boldsymbol{R}_t\}_{t\in\mathbb{Z}}$ is also $\beta$-mixing with $\beta(n) \le n^{-1-\epsilon}$. Thus, by Lemma 7.2, we have

$$\frac{\sqrt{T}A_1}{\sigma} = \frac{\sqrt{T}(U_T - \theta)}{\sigma} \xrightarrow{\mathsf{d}} Z, \tag{7.41}$$

where $Z \sim N(0, 1)$ is a standard Gaussian random variable. By Slutsky's theorem, combining the above equation with (7.38) leads to the desired result. □

## 7.3 Proof of Theorem 3.2

*Proof.* Similar to the proof of Theorem 3.1, we can expand $\mathbf{w}^\mathsf{T}(\widehat{\mathbf{\Sigma}}^* - \mathbf{\Sigma})\mathbf{w}$ by

$$\mathbf{w}^\mathsf{T}(\widehat{\mathbf{\Sigma}}^*-\mathbf{\Sigma})\mathbf{w}=\underbrace{\boldsymbol{a}^\mathsf{T}\Big\{\cos(\frac{\pi}{2}\mathbf{T})\circ\frac{\pi}{2}(\widehat{\mathbf{T}}^*-\mathbf{T})\Big\}\boldsymbol{a}}_{A_1^*}+\underbrace{\boldsymbol{a}^\mathsf{T}\Big\{-\frac{1}{2}[\sin(\theta_{jk})]\circ(\frac{\pi}{2})^2(\widehat{\mathbf{T}}^*-\mathbf{T})\circ(\widehat{\mathbf{T}}^*-\mathbf{T})\Big\}\boldsymbol{a}}_{A_2^*}. \quad (7.42)$$

Let $\widehat{R}^* := \mathbf{w}^\mathsf{T}(\widehat{\mathbf{\Sigma}}^* - \mathbf{\Sigma})\mathbf{w}$ and rewrite $A_1^*$ as

$$A_1^* = \underbrace{\frac{2}{T(T-1)}\sum_{t<t'} g(\boldsymbol{R}_t^*, \boldsymbol{R}_{t'}^*)}_{U_T^*} - \boldsymbol{a}^\mathsf{T}\Big\{\cos(\frac{\pi}{2}\mathbf{T})\circ\frac{\pi}{2}\mathbf{T}\Big\}\boldsymbol{a}.$$

Remind that $g(\cdot)$ is a bounded kernel function and Assumption **(A3)** implies that the process $\{\boldsymbol{R}_t\}_{t\in\mathbb{Z}}$ is $\beta$-mixing with $\beta(n) \le n^{-1-\epsilon}$. By Lemma 7.4 and Assumption **(A2)**, we then have

$$\Big|\mathrm{Var}^*(\sqrt{T}U_T^*) - \mathrm{Var}(\sqrt{T}U_T)\Big| = o_P(\sigma^2),$$

where $U_T$ is defined in Equation (7.40). Moreover, by (7.41), we have $\mathrm{Var}(\sqrt{T}U_T) = \sigma^2\{1 + o(1)\}$. Thus, we have

$$\mathrm{Var}^*(\sqrt{T}A_1^*) = \mathrm{Var}^*(\sqrt{T}U_T^*) = \sigma^2\{1 + o_P(1)\}. \quad (7.43)$$

Next, we focus on the asymptotics of $\mathrm{Var}^*(\sqrt{T}A_2^*)$. Noting that by (7.39), we have

$$A_2^* \le \frac{\pi^2}{4}\|\boldsymbol{a}\|_1^2\|\widehat{\mathbf{T}}^* - \mathbf{T}\|_{\max}^2.$$

By the circular block bootstrap procedure, the process $\{\boldsymbol{R}_t^*\}_{t\in\mathbb{Z}}$ is still a $\phi$-mixing process with mixing coefficient $\widetilde{\phi}(n) \le n^{-(1+\epsilon)(1-\epsilon_0)} = O(n^{-1-\epsilon_2})$ for some $\epsilon_2 > 0$ as long as $\epsilon > \epsilon_0/(1-\epsilon_0)$. Thus, by Lemma 7.6, we have $\|\widehat{\mathbf{T}}^* - \mathbf{T}\|_{\max} = O_P(\sqrt{\log d/T})$. Thus, we have $A_2^* = O_P(\log d/T)$ and accordingly

$$\mathrm{Var}^*(\sqrt{T}A_2^*) \le T\mathbb{E}^*(A_2^{*^2}) = O_P\Big\{\frac{(\log d)^2}{T}\Big\} = o_P(\sigma^2), \quad (7.44)$$

where $\mathbb{E}^*$ is the bootstrap expectation conditional on $\{\boldsymbol{R}_t\}_{t=1}^T$. Combining Equations (7.42), (7.43), and (7.44), we have

$$\begin{aligned}\mathrm{Var}^*(\sqrt{T}\widehat{R}^*) &= \mathrm{Var}^*\Big\{\sqrt{T}(A_1^* + A_2^*)\Big\} = \mathrm{Var}^*(\sqrt{T}A_1^*) + \mathrm{Var}^*(\sqrt{T}A_2^*) + 2\,\mathrm{Cov}(\sqrt{T}A_1^*, \sqrt{T}A_2^*)\\ &\le \mathrm{Var}^*(\sqrt{T}A_1^*) + \mathrm{Var}^*(\sqrt{T}A_2^*) + 2\sqrt{\mathrm{Var}^*(\sqrt{T}A_1^*)}\sqrt{\mathrm{Var}^*(\sqrt{T}A_2^*)}\\ &= \sigma^2\{1 + o_P(1)\}. \quad (7.45)\end{aligned}$$

On the other hand, we also have

$$\begin{aligned}\mathrm{Var}^*(\sqrt{T}\widehat{R}^*) &\geq \mathrm{Var}^*(\sqrt{T}A_1^*) + \mathrm{Var}^*(\sqrt{T}A_2^*) - 2\sqrt{\mathrm{Var}^*(\sqrt{T}A_1^*)}\sqrt{\mathrm{Var}^*(\sqrt{T}A_2^*)}\\ &= \sigma^2\{1 + o_P(1)\}. \end{aligned} \tag{7.46}$$

Combining (7.45) and (7.46) completes the proof. $\square$

## 7.4 Proof of Theorem 3.5

*Proof.* Denote $\boldsymbol{a} := \mathbf{D}\mathbf{w}$. We can write

$$\mathbf{w}^\intercal(\widehat{\boldsymbol{\Sigma}}^s - \boldsymbol{\Sigma})\mathbf{w} = \underbrace{\boldsymbol{a}^\intercal\Big\{\sin(\frac{\pi}{2}\widehat{\mathbf{T}}^s) - \sin(\frac{\pi}{2}\mathbf{T})\Big\}\boldsymbol{a}}_{B_1} + \underbrace{\Big\{\mathbf{w}^\intercal\widehat{\mathbf{D}}\sin(\frac{\pi}{2}\widehat{\mathbf{T}}^s)\widehat{\mathbf{D}}\mathbf{w} - \mathbf{w}^\intercal\mathbf{D}\sin(\frac{\pi}{2}\widehat{\mathbf{T}}^s)\mathbf{D}\mathbf{w}\Big\}}_{B_2}. \tag{7.47}$$

By the same arguments as in the proof of Theorem 3.1, we have

$$\frac{\sqrt{T_s}B_1}{\sigma} \xrightarrow{\mathsf{d}} Z, \tag{7.48}$$

where $Z \sim N(0,1)$ is a Gaussian random variable. It remains to show that $B_2$ is ignorable asymptotically. Using (7.39), we have

$$\begin{aligned}|B_2| &\leq \Big|\mathbf{w}^\intercal\widehat{\mathbf{D}}\sin(\frac{\pi}{2}\widehat{\mathbf{T}}^s)(\widehat{\mathbf{D}} - \mathbf{D})\mathbf{w}\Big| + \Big|\mathbf{w}^\intercal(\widehat{\mathbf{D}} - \mathbf{D})\sin(\frac{\pi}{2}\widehat{\mathbf{T}}^s)\mathbf{D}\mathbf{w}\Big|\\ &\leq \|\sin\frac{\pi}{2}\widehat{\mathbf{T}}\|_{\max}\|(\widehat{\mathbf{D}} - \mathbf{D})\mathbf{w}\|_1(\|\widehat{\mathbf{D}}\mathbf{w}\|_1 + \|\mathbf{D}\mathbf{w}\|_1)\\ &\leq \|\widehat{\mathbf{D}} - \mathbf{D}\|_{\max}(\|\widehat{\mathbf{D}}\|_{\max} + \|\mathbf{D}\|_{\max}).\end{aligned}$$

Using Lemma 7.9 and Assumption **(A5)**, we have $|B_2| = O_P(\sqrt{\log d/T}) = o_P(\sigma/\sqrt{T_s})$. Together with (7.47) and (7.48), using Slutsky's theorem, we have the desired result. $\square$

## 7.5 Proof of Corollary 3.1

*Proof.* By (3.5), we have $\mathbb{P}(|\widehat{w}_j/w_j - 1| > t) \leq \exp(-CTt^2)$. Thus, we further have

$$\mathbb{P}(\max_j|\widehat{w}_j/w_j - 1| > t) \leq d\mathbb{P}(|\widehat{w}_j/w_j - 1| > t) \leq \exp(\log d - CTt^2).$$

To simplify the rate of convergence, setting $t = \sqrt{(3\log d)/(CT)}$, we have

$$\mathbb{P}\Big(\max_j|\widehat{w}_j/w_j - 1| > \sqrt{3\log d/(CT)}\Big) \leq 1/d.$$

Thus, as $(T,d)$ go to infinity, we have $\max_j |\widehat{w}_j/w_j - 1| = O_P(\sqrt{\log d/T})$. This gives us an upper bound of the convergence rate of $\|\widehat{\mathbf{w}} - \mathbf{w}\|_1$:

$$\|\widehat{\mathbf{w}}-\mathbf{w}\|_1=\sum_{j=1}^d |\widehat{w}_j-w_j|=\sum_{j=1}^d |w_j|\cdot\Big|\frac{\widehat{w}_j}{w_j}-1\Big|\le\|\mathbf{w}\|_1\cdot\max_j\Big|\frac{\widehat{w}_j}{w_j}-1\Big| = O_P\Big(\sqrt{\frac{\log d}{T}}\Big). \quad (7.49)$$

Similar as in (7.47), we can decompose $\widehat{\mathbf{w}}^\mathsf{T}\widehat{\mathbf{\Sigma}}^s\widehat{\mathbf{w}} - \mathbf{w}^\mathsf{T}\mathbf{\Sigma}\mathbf{w}$ into

$$\widehat{\mathbf{w}}^\mathsf{T}\widehat{\mathbf{\Sigma}}^s\widehat{\mathbf{w}} - \mathbf{w}^\mathsf{T}\mathbf{\Sigma}\mathbf{w} = B_1 + \underbrace{\widehat{\mathbf{w}}^\mathsf{T}\widehat{\mathbf{D}}\sin(\frac{\pi}{2}\widehat{\mathbf{T}}^s)\widehat{\mathbf{D}}\widehat{\mathbf{w}} - \mathbf{w}^\mathsf{T}\mathbf{D}\sin(\frac{\pi}{2}\widehat{\mathbf{T}}^s)\mathbf{D}\mathbf{w}}_{B_3}, \quad (7.50)$$

where $B_1$ is defined in (7.47). As in the proof of Theorem 3.5, we still have (7.48). Regarding $B_3$, we have $|B_3| \le \|\widehat{\mathbf{D}}\widehat{\mathbf{w}} - \mathbf{D}\mathbf{w}\|_1\|\widehat{\mathbf{D}}\widehat{\mathbf{w}} + \mathbf{D}\mathbf{w}\|_1$. Using the triangle inequality, we have

$$\begin{aligned}|B_3| &\le \Big(\|\widehat{\mathbf{D}}(\widehat{\mathbf{w}} - \mathbf{w})\|_1 + \|(\widehat{\mathbf{D}} - \mathbf{D})\mathbf{w}\|_1\Big)\Big(\|\widehat{\mathbf{D}}\|_1\|\widehat{\mathbf{w}}\|_1 + \|\mathbf{D}\|_1\|\mathbf{w}\|_1\Big)\\ &\le \Big(\|\widehat{\mathbf{D}}\|_{\max}\|\widehat{\mathbf{w}} - \mathbf{w}\|_1 + \|\widehat{\mathbf{D}} - \mathbf{D}\|_{\max}\Big)\Big(\|\widehat{\mathbf{D}}\|_{\max}\|\widehat{\mathbf{w}}\|_1 + \|\mathbf{D}\|_{\max}\Big).\end{aligned}$$

Using (7.49) and Lemma 7.9, we can conclude $|B_3| = O_P(\sqrt{\log d/T})$. Plugging it into (7.50) and using the Slutsky's theorem, we have the desired result. □

## 7.6 Proof of Theorem 3.6

*Proof.* Let $\widehat{\mathbf{K}}^{s^*} = \sin(\pi\widehat{\mathbf{T}}^{s^*}/2)$ and $\mathbf{K} = \sin(\pi\mathbf{T}/2)$. We can decompose $\widehat{R}_s^* := \mathbf{w}^\mathsf{T}\widehat{\mathbf{D}}\widehat{\mathbf{K}}^{s^*}\widehat{\mathbf{D}}\mathbf{w}$ into two parts:

$$\widehat{R}_s^* = \underbrace{\mathbf{w}^\mathsf{T}\mathbf{D}\widehat{\mathbf{K}}^{s^*}\mathbf{D}\mathbf{w} - \mathbf{w}^\mathsf{T}\mathbf{D}\mathbf{K}\mathbf{D}\mathbf{w}}_{B_1^*} + \underbrace{\mathbf{w}^\mathsf{T}\widehat{\mathbf{D}}\widehat{\mathbf{K}}^{s^*}\widehat{\mathbf{D}}\mathbf{w} - \mathbf{w}^\mathsf{T}\mathbf{D}\widehat{\mathbf{K}}^{s^*}\mathbf{D}\mathbf{w}}_{B_2^*}. \quad (7.51)$$

By similar arguments as in the proof of Theorem 3.2, we have

$$\mathrm{Var}^*(\sqrt{T_s}B_1^*) = \sigma^2\{1 + o_P(1)\}. \quad (7.52)$$

Next, we show that $\mathrm{Var}^*(\sqrt{T_s}B_2^*) = o_P(\sigma^2)$. We can upper bound $\mathrm{Var}^*(B_2^*)$ by

$$\begin{aligned}\mathrm{Var}^*(B_2^*) &= \mathrm{Var}^*\Big\{\mathbf{w}^\mathsf{T}\widehat{\mathbf{D}}\widehat{\mathbf{K}}^{s^*}(\widehat{\mathbf{D}} - \mathbf{D})\mathbf{w} + \mathbf{w}^\mathsf{T}(\widehat{\mathbf{D}} - \mathbf{D})\widehat{\mathbf{K}}^{s^*}\mathbf{D}\mathbf{w}\Big\}\\ &\le \mathrm{Var}^*\Big\{\mathbf{w}^\mathsf{T}\widehat{\mathbf{D}}\widehat{\mathbf{K}}^{s^*}(\widehat{\mathbf{D}} - \mathbf{D})\mathbf{w}\Big\} + \mathrm{Var}^*\Big\{\mathbf{w}^\mathsf{T}(\widehat{\mathbf{D}} - \mathbf{D})\widehat{\mathbf{K}}^{s^*}\mathbf{D}\mathbf{w}\Big\}\\ &\quad + 2\sqrt{\mathrm{Var}^*\Big\{\mathbf{w}^\mathsf{T}\widehat{\mathbf{D}}\widehat{\mathbf{K}}^{s^*}(\widehat{\mathbf{D}} - \mathbf{D})\mathbf{w}\Big\}}\sqrt{\mathrm{Var}^*\Big\{\mathbf{w}^\mathsf{T}(\widehat{\mathbf{D}} - \mathbf{D})\widehat{\mathbf{K}}^{s^*}\mathbf{D}\mathbf{w}\Big\}}. \end{aligned} \quad (7.53)$$

For any random matrix $\mathbf{X} := (\boldsymbol{R}_1, \dots, \boldsymbol{R}_m)^\mathsf{T} \in \mathbb{R}^{m\times n}$ and fixed vectors $\mathbf{v}_1 \in \mathbb{R}^m$, $\mathbf{v}_2 \in \mathbb{R}^n$, let $\mathbf{V}$ be a matrix with $(j,k)$ entry $\mathbf{v}_2^\mathsf{T} \operatorname{Cov}(\boldsymbol{R}_j, \boldsymbol{R}_k)\mathbf{v}_2$. It is easy to verify that

$$\begin{aligned}\operatorname{Var}(\mathbf{v}_1^\mathsf{T}\mathbf{X}\mathbf{v}_2) =& \mathbf{v}_1^\mathsf{T} \operatorname{Var}(\mathbf{X}\mathbf{v}_2)\mathbf{v}_1 = \mathbf{v}_1^\mathsf{T}\mathbf{V}\mathbf{v}_1 \le \|\mathbf{v}_1\|_1^2 \max_{jk}\left|\mathbf{v}_2^\mathsf{T} \operatorname{Cov}(\boldsymbol{R}_j, \boldsymbol{R}_k)\mathbf{v}_2\right| \\ \le& \|\mathbf{v}_1\|_1^2\|\mathbf{v}_2\|_1^2 \max_{j_1,k_1,j_2,k_2} |\operatorname{Cov}(\boldsymbol{R}_{j1,k1}, \boldsymbol{R}_{j2,k2})|. \end{aligned} \tag{7.54}$$

Now writing $\mathbf{v}_1 = \widehat{\mathbf{D}}\mathbf{w}$, $\mathbf{v}_2 = (\widehat{\mathbf{D}} - \mathbf{D})\mathbf{w}$, and $\mathbf{X} = \widehat{\mathbf{K}}^{s^*}$, we have

$$\begin{aligned}\operatorname{Var}^*\left\{\mathbf{w}^\mathsf{T}\widehat{\mathbf{D}}\widehat{\mathbf{K}}^{s^*}(\widehat{\mathbf{D}} - \mathbf{D})\mathbf{w}\right\} \le& \|\widehat{\mathbf{D}}\mathbf{w}\|_1^2\|(\widehat{\mathbf{D}} - \mathbf{D})\mathbf{w}\|_1^2 \max_{j_1,k_1,j_2,k_2} |\overset{*}{\operatorname{Cov}}(\widehat{\tau}^{s^*}_{j_1,k_1}, \widehat{\tau}^{s^*}_{j_2,k_2})| \\ \le& \|\mathbf{w}\|_1^4\|\widehat{\mathbf{D}}\|^2_{\max}\|\widehat{\mathbf{D}} - \mathbf{D}\|^2_{\max} = \|\widehat{\mathbf{D}}\|^2_{\max}\|\widehat{\mathbf{D}} - \mathbf{D}\|^2_{\max}. \end{aligned} \tag{7.55}$$

Note that $\widehat{\mathbf{D}}$ only depends on $\{\boldsymbol{R}_t\}_{t=1}^T$ and is thus fixed under $\operatorname{Var}^*(\cdot)$. Using Lemma 7.9 and (7.55), we have

$$\operatorname{Var}^*\left\{\sqrt{T_s}\mathbf{w}^\mathsf{T}\widehat{\mathbf{D}}\widehat{\mathbf{K}}^{s^*}(\widehat{\mathbf{D}} - \mathbf{D})\mathbf{w}\right\} = O_P\Big(T_s\frac{\log d}{T}\Big) = O_P\Big(\frac{\log d}{T^\delta}\Big) = o_P(\sigma^2). \tag{7.56}$$

Similarly, we also have

$$\operatorname{Var}^*\left\{\sqrt{T_s}\mathbf{w}^\mathsf{T}(\widehat{\mathbf{D}} - \mathbf{D})\widehat{\mathbf{K}}^{s^*}\mathbf{D}\mathbf{w}\right\} = o_P(\sigma^2). \tag{7.57}$$

Combining (7.53), (7.56), and (7.57), we have

$$\operatorname{Var}^*(\sqrt{T_s}B_2^*) = o_P(\sigma^2). \tag{7.58}$$

By (7.51), we have

$$\operatorname{Var}^*(\sqrt{T_s}\widehat{R}^*) \ge \operatorname{Var}^*(\sqrt{T_s}B_1^*) + \operatorname{Var}^*(\sqrt{T_s}B_2^*) - 2\sqrt{\operatorname{Var}^*(\sqrt{T_s}B_1^*)}\sqrt{\operatorname{Var}^*(\sqrt{T_s}B_2^*)},$$

and similarly

$$\operatorname{Var}^*(\sqrt{T_s}\widehat{R}^*) \le \operatorname{Var}^*(\sqrt{T_s}B_1^*) + \operatorname{Var}^*(\sqrt{T_s}B_2^*) + 2\sqrt{\operatorname{Var}^*(\sqrt{T_s}B_1^*)}\sqrt{\operatorname{Var}^*(\sqrt{T_s}B_2^*)}.$$

Using the above two inequalities with (7.52) and (7.58), we can conclude that $\operatorname{Var}^*(\sqrt{T_s}\widehat{R}^*) = \sigma^2\{1 + o_P(1)\}$. □

## 7.7 Proofs of Theorems 3.3 and 3.4

The proofs of Theorems 3.3 and 3.4 are close to those of Theorems 3.5 and 3.6. The main difference is that now $T_h$ plays the role of $T$, and $T$ plays the role of $T_s$. We accordingly omit the proofs.

# References

Agarwal, A., Negahban, S., and Wainwright, M. J. (2012). Noisy matrix decomposition via convex relaxation: Optimal rates in high dimensions. *The Annals of Statistics*, 40(2):1171–1197.

Bai, J. and Liao, Y. (2012). Efficient estimation of approximate factor models via regularized maximum likelihood. *arXiv preprint arXiv:1209.5911*.

Barndorff-Nielsen, O. E. (2002). Econometric analysis of realized volatility and its use in estimating stochastic volatility models. *Journal of the Royal Statistical Society: Series B (Statistical Methodology)*, 64(2):253–280.

Bickel, P. J. and Levina, E. (2008a). Covariance regularization by thresholding. *The Annals of Statistics*, 36(6):2577–2604.

Bickel, P. J. and Levina, E. (2008b). Regularized estimation of large covariance matrices. *The Annals of Statistics*, 36(1):199–227.

Bollerslev, T. (1986). Generalized autoregressive conditional heteroskedasticity. *Journal of Econometrics*, 31(3):307–327.

Cai, T. T., Zhang, C.-H., and Zhou, H. H. (2010). Optimal rates of convergence for covariance matrix estimation. *The Annals of Statistics*, 38(4):2118–2144.

Cai, T. T. and Zhou, H. H. (2012). Optimal rates of convergence for sparse covariance matrix estimation. *The Annals of Statistics*, 40(5):2389–2420.

Chamberlain, G. (1983). A characterization of the distributions that imply mean–variance utility functions. *Journal of Economic Theory*, 29(1):185–201.

Chang, C. and Tsay, R. S. (2010). Estimation of covariance matrix via the sparse Cholesky factor with lasso. *Journal of Statistical Planning and Inference*, 140(12):3858–3873.

Chen, X., Xu, M., and Wu, W. (2013). Covariance and precision matrix estimation for high-dimensional time series. *The Annals of Statistics*, 41(6):2994–3021.

Cont, R. (2001). Empirical properties of asset returns: stylized facts and statistical issues. *Quantatitive Finance*, 1(2):223–236.

Doukhan, P. and Louhichi, S. (1999). A new weak dependence condition and applications to moment inequalities. *Stochastic Processes and their Applications*, 84(2):313–342.

Doukhan, P. and Neumann, M. H. (2007). Probability and moment inequalities for sums of weakly dependent random variables, with applications. *Stochastic Processes and their Applications*, 117(7):878–903.

Fan, J., Fan, Y., and Lv, J. (2008). High dimensional covariance matrix estimation using a factor model. *Journal of Econometrics*, 147(1):186–197.

Fan, J., Han, F., and Liu, H. (2014). PAGE: Robust pattern guided estimation of large covariance matrix. Technical report, Princeton University.

Fan, J., Liao, Y., and Mincheva, M. (2011). High dimensional covariance matrix estimation in approximate factor models. *The Annals of Statistics*, 39(6):3320–3356.

Fan, J., Liao, Y., and Mincheva, M. (2013). Large covariance estimation by thresholding principal orthogonal complements. *Journal of the Royal Statistical Society: Series B (Statistical Methodology)*, 75(4):603–680.

Fan, J., Liao, Y., and Shi, X. (2015). Risks of large portfolios. *Journal of Econometrics (to appear)*.

Fan, J. and Peng, H. (2004). Nonconcave penalized likelihood with a diverging number of parameters. *The Annals of Statistics*, 32(3):928–961.

Fan, J., Zhang, J., and Yu, K. (2012). Vast portfolio selection with gross-exposure constraints. *Journal of the American Statistical Association*, 107(498):592–606.

Frahm, G. and Jaekel, U. (2007). Tyler's M-estimator, random matrix theory, and generalized elliptical distributions with applications to finance. Technical report, Helmut Schmidt University.

Fryzlewicz, P. (2013). High-dimensional volatility matrix estimation via wavelets and thresholding. *Biometrika*, 100(4):921–938.

Gómez, K. and Gallón, S. (2011). Comparison among high dimensional covariance matrix estimation methods. *Revista Colombiana de Estadística*, 34(3):567–588.

Greenshtein, E. and Ritov, Y. (2004). Persistence in high-dimensional linear predictor selection and the virtue of overparametrization. *Bernoulli*, 10(6):971–988.

Hamada, M. and Valdez, E. (2004). *CAPM and option pricing with elliptical distributions*. School of Finance and Economics, University of Technology, Sydney.

Han, F. and Liu, H. (2013a). Optimal rates of convergence for latent generalized correlation matrix estimation in transelliptical distribution. *arXiv preprint arXiv:1305.6916*.

Han, F. and Liu, H. (2013b). Principal component analysis on non-Gaussian dependent data. In *Proceedings of the 30th International Conference on Machine Learning*, pages 240–248.

Han, F. and Liu, H. (2013c). Transition matrix estimation in high dimensional time series. In *Proceedings of the 30th International Conference on Machine Learning*, pages 172–180.

Han, F. and Liu, H. (2014a). Distribution-free tests of independence with applications to testing more structures. *arXiv preprint arXiv:1410.4179*.

Han, F. and Liu, H. (2014b). Scale-invariant sparse PCA on high-dimensional meta-elliptical data. *Journal of the American Statistical Association*, 109(505):275–287.

Han, F., Lu, J., and Liu, H. (2014). Robust scatter matrix estimation for high dimensional distributions with heavy tails. Technical report, Princeton University.

Higham, N. J. (2002). Computing the nearest correlation matrix—a problem from finance. *IMA Journal of Numerical Analysis*, 22(3):329–343.

Hsu, D., Kakade, S. M., and Zhang, T. (2011). Robust matrix decomposition with sparse corruptions. *IEEE Transactions on Information Theory*, 57(11):7221–7234.

Jagannathan, R. and Ma, T. (2003). Risk reduction in large portfolios: Why imposing the wrong constraints helps. *The Journal of Finance*, 58(4):1651–1684.

Kontorovich, L. A., Ramanan, K., et al. (2008). Concentration inequalities for dependent random variables via the martingale method. *The Annals of Probability*, 36(6):2126–2158.

Lai, T. L., Xing, H., and Chen, Z. (2011). Mean-variance portfolio optimization when means and covariances are unknown. *The Annals of Applied Statistics*, 5(2A):798–823.

Ledoit, O. and Wolf, M. (2003). Improved estimation of the covariance matrix of stock returns with an application to portfolio selection. *Journal of Empirical Finance*, 10(5):603–621.

Lindskog, F., McNeil, A., and Schmock, U. (2003). Kendall's tau for elliptical distributions. *Credit Risk: Measurement, Evaluation and Management*, pages 149–156.

Loh, P.-L. and Wainwright, M. J. (2012). High-dimensional regression with noisy and missing data: Provable guarantees with nonconvexity. *The Annals of Statistics*, 40(3):1637–1664.

Mitra, R. and Zhang, C.-H. (2014). Multivariate analysis of nonparametric estimates of large correlation matrices. *arXiv preprint arXiv:1403.6195*.

Mohri, M. and Rostamizadeh, A. (2010). Stability bounds for stationary $\phi$-mixing and $\beta$-mixing processes. *The Journal of Machine Learning Research*, 11:789–814.

Owen, J. and Rabinovitch, R. (1983). On the class of elliptical distributions and their applications to the theory of portfolio choice. *The Journal of Finance*, 38(3):745–752.

Pan, J. and Yao, Q. (2008). Modelling multiple time series via common factors. *Biometrika*, 95(2):365–379.

Pesaran, M. H. and Zaffaroni, P. (2008). Optimal asset allocation with factor models for large portfolios. Technical report, CESifo working paper.

Politis, D. N. and Romano, J. P. (1992). A circular block-resampling procedure for stationary data. In *Exploring the Limits of Bootstrap*, pages 263–270. John Wiley, New York.

Qiu, H., Han, F., Liu, H., and Caffo, B. (2014). Robust portfolio optimization under high dimensional heavy-tailed time series. Technical report, Johns Hopkins University.

Shao, Q.-M. and Yu, H. (1993). Bootstrapping the sample means for stationary mixing sequences. *Stochastic Processes and their Applications*, 48(1):175–190.

Wegkamp, M. and Zhao, Y. (2013). Adaptive estimation of the copula correlation matrix for semiparametric elliptical copulas. *arXiv preprint arXiv:1305.6526*.

Xiao, H. and Wu, W. B. (2012). Covariance matrix estimation for stationary time series. *The Annals of Statistics*, 40(1):466–493.

Yoshihara, K.-I. (1976). Limiting behavior of U-statistics for stationary, absolutely regular processes. *Probability Theory and Related Fields*, 35(3):237–252.

Zhang, L., Mykland, P. A., and Aït-Sahalia, Y. (2005). A tale of two time scales. *Journal of the American Statistical Association*, 100(472):1394–1411.